\documentclass[a4paper, 10pt]{amsart}
\usepackage{amsfonts,amsmath,amssymb}
\usepackage{mathrsfs,mathtools}
\usepackage{enumerate}
\usepackage{xspace}
\usepackage{./mydef}

\usepackage{hyperref}
\usepackage{esint}
\usepackage{graphicx}
\usepackage{stmaryrd}
\DeclareMathAlphabet{\mathpzc}{OT1}{pzc}{m}{it}

\newif\ifdraft
\drafttrue
\ifdraft
\usepackage[usenames,dvipsnames]{color}
\fi

\newtheorem{theorem}{Theorem}[section]
\newtheorem{lemma}[theorem]{Lemma}
\newtheorem{proposition}[theorem]{Proposition}

\newtheorem*{proposition*}{Proposition}
\newtheorem*{theorem*}{Theorem}
\theoremstyle{remark}
\newtheorem{remark}[theorem]{Remark}
\theoremstyle{definition}
\newtheorem{definition}[theorem]{Definition}
\newtheorem{example}[theorem]{Example}

\numberwithin{equation}{section}

\newcommand{\dist}{{\textrm{dist}}}
\newcommand{\calG}{{\mathcal G}}

\DeclareMathOperator*{\supp}{supp}
\newcommand{\ve}{\mathscr{V}}
\newcommand{\calI}{\mathcal{I}}

\newcommand{\rn}{{\mathbb{R}^d}}

\newcommand{\pp}{\partial}
\newcommand{\eps}{\varepsilon}
\newcommand{\contour}{{\mathscr{C}}}
\newcommand{\gain}{\alpha}
\newcommand{\SZ}{\Pi_\T}
\newcommand{\wero}{\texttt{w}}

\newcommand{\othervar}{x'}

\begin{document}
%%%%%%%%%%%%%%%%%%%%%%%%%%%%%%%%%%%%%%%%%%%%%%%%%%%%%%%%%%%%%%%%%%%%%%%%%%%%%%%%%
\title[Numerical Methods for Fractional Diffusion]{Numerical Methods
  for Fractional Diffusion}

\author[A.~Bonito]{Andrea Bonito}
\address[A.~Bonito]{Department of Mathematics, Texas A\&M University, College Station, TX 77843, USA}
\thanks{AB is supported in part by NSF grant DMS-1254618.}
\email{bonito@math.tamu.edu}
\author[J.P.~Borthagaray]{Juan Pablo~Borthagaray}
\address[J.P.~Borthagaray]{IMAS - CONICET and Departamento de Matem\'atica, FCEyN - Universidad de Buenos Aires,  Buenos Aires, Argentina}
\email{jpbortha@dm.uba.ar}
\thanks{JPB has been partially supported by a CONICET doctoral fellowship}

\author[R.H.~Nochetto]{Ricardo H.~Nochetto}
\address[R.H.~Nochetto]{Department of Mathematics and Institute for Physical Science and Technology,
University of Maryland, College Park, MD 20742, USA}
\email{rhn@math.umd.edu}
\thanks{RHN has been supported in part by NSF grant DMS-1411808.}

\author[E.~Ot\'arola]{Enrique Ot\'arola} %% El pelaito vaca
\address[E.~Ot\'arola]{Departamento de Matem\'atica, Universidad T\'ecnica Federico Santa Mar\'ia, Valpara\'iso, Chile}
\email{enrique.otarola@usm.cl}
\thanks{EO has been supported in part by CONICYT through project FONDECYT 3160201.}

\author[A.J.~Salgado]{Abner J.~Salgado}
\address[A.J.~Salgado]{Department of Mathematics, University of Tennessee, Knoxville, TN 37996, USA}
\email{asalgad1@utk.edu}
\thanks{AJS is supported by NSF grant DMS-1418784.}

%%%%%%%%%%%%%%%%%%%%%%%%%%%%%%%%%%%%%%%%%%%%%%%%%%%%%%%%%%%%%%%%%%%%%%%%%%%%%%%%%%%

\pagestyle{myheadings}
\thispagestyle{plain}
\markboth{\textsc{A. Bonito, J.P.~Borthagaray, R.H.~Nochetto,
    E.~Ot\'arola, A.J.~Salgado}}{Three numerical methods for fractional diffusion}

\date{Draft version of \today.}

%%%%%%%%%%%%%%%%%%%%%%%%%%%%%%%%%%%%%%%%%%%%%%%%%%%%%%%%%%%%%%%%%%%%%%%%%%%%%%%%%%%

\begin{abstract}
We present three schemes for the numerical approximation of
fractional diffusion, which build on different definitions of such 
a non-local process.
The first method is a PDE approach that applies to the spectral
definition and exploits the extension to one higher dimension.
The second method is the integral formulation and deals with singular
non-integrable kernels. The third method is a discretization of the Dunford-Taylor
formula. We discuss pros and cons of each method, error estimates, and
document their performance with a few numerical experiments.
\end{abstract}

\maketitle

%%%%%%%%%%%%%%%%%%%%%%%%%%%%%%%%%%%%%%%%%%%%%%%%%%%%%%%%%%%%%%%%%%%%%%%%%%%%%%%%%%%%
\section{Introduction}\label{S:Introduction}
%%%%%%%%%%%%%%%%%%%%%%%%%%%%%%%%%%%%%%%%%%%%%%%%%%%%%%%%%%%%%%%%%%%%%%%%%%%%%%%%%%%%

%\AJSc{Henceforth the dimension shall be denoted by $d$}

Diffusion is the tendency of a substance to evenly spread into an available space, and is one of the most common physical processes. The classical models of diffusion, which usually start from the assumption of Brownian motion \cite{Abe2005403,Albert}, lead to well known models and even better studied equations. However, in recent times, it has become evident that many of the assumptions that lead to these models are not always satisfactory or even realistic at all. For this reason, new models of diffusion have been introduced. These, as a rule, are not based on the postulate that the underlying stochastic process is given by Brownian motion, so that the diffusion is regarded as anomalous \cite{MR2090004}. The evidence of anomalous diffusion processes has been reported in physical and social environments, and corresponding models have been proposed in various areas such as electromagnetic fluids \cite{McCN:81}, 
%biology \cite{tetas:)}, 
ground-water solute transport \cite{CG:93}, biology \cite{biology}, finance \cite{Carr.Geman.ea2002}, human travel \cite{Brockmann} and predator search patterns \cite{Sims}.

Of the many possible models of anomalous diffusion, we shall be interested here in so-called fractional diffusion, which is a nonlocal process: to evaluate fractional diffusion at a spatial point, information involving all spatial points is needed. Recently, the analysis of such operators has received a tremendous attention: fractional diffusion has been one of the most studied topics in the past decade \cite{CS:07,MR2680400,MR2270163}.
% The study of nonlocal operators has been an active area of research in different branches of mathematics, and 
% These operators have been employed to model situations in which different length scales are involved. 

The main goal of this work is to review different techniques to approximate the solution of problems involving fractional diffusion. To make matters precise, we will consider the fractional powers of the Dirichlet Laplace operator $\Laps$, which we will simply call the fractional Laplacian. Given $0<s<1$, a bounded Lipschitz domain $\Omega\subset\R^d$, and a function $f:\Omega\to\mathbb{R}$, we shall be concerned with finding $u: \Omega \to \R$ such that
\begin{equation}
\label{fl=f_bdddom}
    \Laps u = f \quad \text{in } \Omega,
\end{equation}
with vanishing Dirichlet boundary conditions (understood in a suitable
sense). We must immediately comment that the efficient approximation
of solutions to \eqref{fl=f_bdddom} carries two essential
difficulties. The first and most important one is that $\Laps$ is \emph{non-local}. The second feature is the lack of boundary regularity, which leads to reduced rates of convergence.

%We consider two different definitions of fractional Laplacian $\Laps u$;
%their connection is studied in \cite{CaffarelliStinga:15}.
In what follows we review different definitions for the fractional Laplacian. For functions defined over $\rn$, there is a natural way to define the fractional Laplacian as a pseudo-differential operator with symbol $|\xi|^{2s}$; namely, given a function $w$ in the Schwartz class $\mathscr{S}$, set
\begin{equation}\label{eq:Fourier}
(-\Delta)^s w := \mathscr{F}^{-1} \left( |\xi|^{2s} \mathscr{F} w \right) ,
\end{equation}
where $\mathscr{F}$ denotes the Fourier transform. The fractional Laplacian can be equivalently defined by means of the following point-wise formula (see \cite[Section 1.1]{Landkof} and \cite[Proposition 3.3]{Hitchhikers})
\begin{equation}
(-\Delta)^s w(x) = C(d,s) \mbox{ p.v.} \int_\rn \frac{w(x)-w(\othervar)}{|x-\othervar|^{d+2s}} \, \diff \othervar, \quad C(d,s) = \frac{2^{2s} s \Gamma(s+\frac{d}{2})}{\pi^{d/2}\Gamma(1-s)},
\label{eq:frac_lap}
\end{equation}
where \mbox{p.v.} stands for the Cauchy principal value and $C(d,s)$ is a normalization constant chosen so that definitions \eqref{eq:Fourier} and \eqref{eq:frac_lap} coincide. This clearly displays the non-local structure of $\Laps w$. We remark that, in the theory of stochastic processes, expression \eqref{eq:frac_lap} appears as the infinitesimal generator of a $2s$-stable L\'evy process \cite{Bertoin}.

If $\Omega$ is a bounded domain, we consider two possible definitions of the fractional Laplacian. 
For $u:\Omega \to \R$, we first extend $u$ by zero outside $\Omega$ and next use definition \eqref{eq:frac_lap}. This gives the following reinterpretation of \eqref{fl=f_bdddom}:
\begin{equation}
  (-\Delta)^s \tilde u = f \quad  \mbox{in }\Omega,
  \qquad
   \tilde u = 0 \quad \mbox{in }\Omega^c = \rn \setminus \Omega ,
\label{eq:Dir_frac}
\end{equation}
where the operator $\Laps$ is understood as in \eqref{eq:frac_lap}
and $\tilde w$ is the extension by zero to $\rn$ of a function $w: \Omega \to \R$ in $L^2(\Omega)$.
This definition maintains the probabilistic interpretation of the fractional Laplacian defined over $\rn$, that is, as the generator of a random walk in $\Omega$ with arbitrarily long jumps, where particles are killed upon reaching $\Omega^c$; see \cite[Chapter 2]{BV}. The operator in \eqref{eq:frac_lap} is well defined for smooth, compactly supported functions. Consequently, \eqref{eq:frac_lap} can be extended by density to $\Hs$ which, for $s \in [0,3/2)$, is defined by
\begin{equation}\label{e:tildeHs}
    \Hs := \left\lbrace w|_{\Omega}: w \in H^s(\rn) , \ w|_{\mathbb{R}^d \setminus \Omega} = 0 \right\}.
\end{equation}
When $\partial\Omega$ is Lipschitz  this space is equivalent to $\Hs = [L^2(\Omega),H_0^1(\Omega)]_s$, the real interpolation between $L^2(\Omega)$ and $H_0^1(\Omega)$ \cite{McLean,Tartar} when $s\in[0,1]$ and to $H^s(\Omega) \cap H^1_0(\Omega)$ when $s \in (1,3/2)$. In what follows, we will denote by $\Hsd$ the dual of $\Hs$ and by $\langle \cdot, \cdot \rangle$ the duality pairing between these two spaces.

The second definition of $(-\Delta)^s$ relies on spectral theory
\cite{BS}. Since $-\Delta: \mathcal{D}(-\Delta) \subset L^2(\Omega) \to
L^2(\Omega)$ is an unbounded, positive and closed operator with dense
domain $\mathcal{D}(-\Delta) = H_0^1(\Omega) \cap H^2(\Omega)$ and its
inverse is compact, there is a countable collection of eigenpairs
$\{\lambda_k, \varphi_k \}_{k \in \mathbb{N}} \subset \R^+ \times
H_0^1(\Omega)$ such that $\{\varphi_k\}_{k \in \mathbb{N}}$ is an
orthonormal basis of $L^2(\Omega)$ as well as an orthogonal basis of $H_0^1(\Omega)$.
Fractional powers of the Dirichlet Laplacian can be thus defined as
\begin{equation}
  \label{def:second_frac}
  (-\Delta)^s w  := \sum_{k=1}^\infty \lambda_k^{s} w_k \varphi_k,  \quad w_k := \int_{\Omega} w \varphi_k \diff x, \quad k \in \mathbb{N},
\end{equation}
for any $w \in C_0^\infty(\Omega)$.
This definition of $(-\Delta)^s$ can also be extended by density to the space $\Hs$. We also remark that, for $s \in [-1,1]$, the space $\Hs$ can equivalently be defined by
\[
  \Hs = \left\{ w = \sum_{k=1}^\infty w_k \varphi_k \ :  \ \sum_{k=1}^\infty \lambda_k^s w_k^2 < \infty \right\} .
\]

These two definitions of the fractional Laplacian, the \emph{integral} one involved in problem \eqref{eq:Dir_frac} and the \emph{spectral} one given as in \eqref{def:second_frac}, do not coincide. In fact, as shown in \cite{MR3246044}, their difference is positive and positivity preserving, see also \cite{CaffarelliStinga:15,MR3233760}. This, in particular, implies that the boundary behavior of the solutions of \eqref{eq:Dir_frac} and \eqref{def:second_frac}
is quite different. According to Grubb \cite{Grubb} the solution $u$
of \eqref{eq:Dir_frac} is of the form
\begin{equation}\label{boundary-grubb}
u(x) \approx \textrm{dist} (x,\partial\Omega)^s + v(x),
\end{equation}
with $v$ smooth; hereafter $\textrm{dist} (x,\partial\Omega)$
indicates the distance from $x\in\Omega$ to $\partial\Omega$.
In contrast, Caffarelli and Stinga \cite{CaffarelliStinga:15} showed that
solutions of \eqref{def:second_frac} behave like
\begin{equation}\label{boundary-CS}
\begin{aligned}
u(x) &\approx \textrm{dist} (x,\partial\Omega)^{2s} + v(x)
\qquad 0 <s <\frac12;
\\
u(x) &\approx \textrm{dist} (x,\partial\Omega)+ v(x)
\qquad \frac12 \le s < 1.
\end{aligned}
\end{equation}
This lack of boundary regularity is responsible for reduced rates of convergence.

The presence of the non-integrable kernel in \eqref{eq:frac_lap} is a
notorious numerical difficulty that has hampered progress in the
multidimensional case $d>1$ until recently. If $\Omega=\rn$,
the Caffarelli-Silvestre \cite{CS:07}
extension converts \eqref{fl=f_bdddom} into the following
Dirichlet-to-Neumann map formulated in the cylinder
$\C=\Omega\times(0,\infty)$
\begin{equation}
\label{alpha_harm_intro}
  \DIV\left( y^\alpha \nabla \ue \right) = 0  \text{ in } \C, \qquad
  \partial_{\nu^\alpha} \ue = d_s f  \text{ on } \rn \times \{0\},
\end{equation}
where $y>0$ is the extended variable, $d_s: = 2^{1-2s}\Gamma(1-s)/\Gamma(s)$
is a positive normalization constant that
depends only on $s$, and the parameter $\alpha$ is defined as $\alpha :=
1-2s \in (-1,1)$. The relation between \eqref{fl=f_bdddom} and
\eqref{alpha_harm_intro} is the following:
\[
u=\ue(\cdot,0).
\]
Cabr\'e and Tan \cite{CT:10} and Stinga and Torrea \cite{ST:10} have shown that a similar extension is valid for the
spectral fractional Laplacian in a bounded domain $\Omega$
provided that a vanishing Dirichlet condition is
appended on the lateral boundary $\partial_L \C = \partial \Omega \times
(0,\infty)$; see also \cite{CDDS:11,BCdPS:12}. Although \eqref{alpha_harm_intro} is a local problem, and
thus amenable to PDE techniques, it
is formulated in one higher dimension and exhibits a singular
character as $y\downarrow0$.

The solution to \eqref{fl=f_bdddom} with either definition \eqref{eq:frac_lap} or \eqref{def:second_frac} of the fractional Laplacian can be represented using Dunford-Taylor integrals \cite{lunardi}. 
% Yoshida reference does not apply
Let us first explain the construction for definition \eqref{def:second_frac}. For $s \in (0,1)$ and $w \in \Hsd$ we have
\begin{equation}
\label{e:dunf-inv}
  (-\Delta)^{-s}w = \frac{1}{2\pi i} \int_{\contour} z^{-s} (z+\Delta)^{-1}w \diff z,
\end{equation}
where $\contour$ is a Jordan curve oriented to have the spectrum of $-\Delta$ to its right, and $z^{-s}$ is defined using the principal value of $\mathrm{Log}(z)$. In addition, since the operator $-\Delta$ is positive, one can continously deform the contour $\contour$ onto the negative real axis (around the branch cut) to obtain the so-called Balakrishnan formula 
\begin{equation}
\label{e:balak_inv}
  (-\Delta)^{-s}w = \frac{\sin(s \pi)}{\pi } \int_{0}^\infty \mu^{-s} (\mu-\Delta)^{-1}w \diff\mu;
\end{equation}
see also \cite[Section IX.11]{MR0239384} and \cite[Section 10.4]{BS} for a different derivation of \eqref{e:balak_inv} using semigroup theory. This formula will be the starting point for our Dunford-Taylor approach for the spectral fractional Laplacian. These considerations, however, cannot be carried out for the the integral fractional Laplacian \eqref{eq:frac_lap} since neither  \eqref{e:dunf-inv} or \eqref{e:balak_inv} are well defined quantities for this operator. Therefore we will, instead, multiply \eqref{eq:Dir_frac} by a test function $w$, integrate over $\rn$ and use Parseval's equality to obtain the following weak formulation of \eqref{eq:Dir_frac}: Given $f \in \Hsd$ find $u \in \Hs$ such that
\begin{equation}\label{e:weak_int_lap}
  \int_{\rn} |\xi|^s \mathscr{F}^{-1}(\tilde u)  |\xi|^s \overline{\mathscr{F}^{-1}(\tilde w)} \diff \xi =\langle  f, w \rangle, \quad \forall w \in \Hs,
\end{equation}
where $\overline{z}$  denotes the complex conjugate of $z\in \mathbb C$.
Using again the Fourier transform and Parseval's equality, the left hand side of the above relation can be equivalently written as (see Theorem~\ref{T:dunf:int:equv})
\begin{equation}
\label{e:rep_int}
  \frac{2\sin(s\pi)}{\pi} \int_0^\infty \mu^{1-2s} \int_{\Omega} \big((-\Delta)(I-\mu^2 \Delta)^{-1}\tilde u(x) \big) w(x) \diff x \diff \mu.
\end{equation}
These ideas will be the starting point of the Dunford-Taylor method for the integral fractional Laplacian \eqref{eq:frac_lap}.

The purpose of this paper is to present and briefly analyze three finite element methods (FEMs) to approximate \eqref{fl=f_bdddom}. 
The first will consider the \emph{spectral} definition \eqref{def:second_frac}, the second will deal with the \emph{integral} definition \eqref{eq:frac_lap}, while the third approach will be able to account for both operators by means of either \eqref{e:balak_inv} or \eqref{e:rep_int}.
We must immediately remark that for special domain geometries, such as when $d=2$ and $\Omega$ is a rectangle, the use of
spectral methods can be quite efficient but we do not elaborate any further as we are interested in techniques that apply to general domains. Our presentation is organized as follows:
In Section \ref{S:SpectralLaplacian}, we present a method that hinges on the extension \eqref{alpha_harm_intro} and uses PDE techniques. The second method deals with the integral formulation and is discussed in Section \ref{S:IntegralMethod}.
Finally, the third method is based on exponentially convergent quadrature approximations of  \eqref{e:balak_inv} for the spectral Laplacian and of \eqref{e:rep_int} for the integral Laplacian. They are discussed in Section \ref{S:DumfordTaylor}. 

As usual, we write $a \lesssim b$ to mean $a \leq Cb$, with a constant
$C$ that neither depends on $a,b$ or the discretization
parameters and might change at each occurrence. Moreover, $a \approx b$
indicates $a \lesssim b$ and $b \lesssim a$.

% !TEX root = three-methods.tex
%%%%%%%%%%%%%%%%%%%%%%%%%%%%%%%%%%%%%%%%%%%%%%%%%%%%%%%%%%%%%%%%%%%%%%%%%%%%%%%%%%%%
\section{The Spectral Fractional Laplacian}\label{S:SpectralLaplacian}
%%%%%%%%%%%%%%%%%%%%%%%%%%%%%%%%%%%%%%%%%%%%%%%%%%%%%%%%%%%%%%%%%%%%%%%%%%%%%%%%%%%%

In this section we deal with the spectral definition of the fractional
Laplacian \eqref{def:second_frac} and, on the basis of
\eqref{alpha_harm_intro}, its discretization via PDE techniques as
originally developed in \cite{NOS}. We must immediately remark that
many of the results of this section and section \ref{S:DumfordTaylor} extend to more general symmetric elliptic operators of the form
$
L w = - \DIV( A \GRAD w ) + c w
$,
with $A \in C^{0,1}(\bar\Omega,\GL(\R^d))$ symmetric and positive definite 
and $0 \leq c \in C^{0,1}(\bar\Omega,\R)$.

Let $\Omega$ be a convex polytopal domain.
Besides the semi-infinite cylinder $\C= \Omega\times (0,\infty)$, we introduce
the truncated cylinder $\C_\Y = \Omega \times (0,\Y)$ with height $\Y$
and its lateral boundary $\partial_L\C_\Y = \partial \Omega \times (0,\Y)$.
Since we deal with objects defined in both $\R^d$ and $\R^{d+1}$, it
is convenient to distinguish the extended variable $y$. For
$\bx\in\R^{d+1}$, we denote
\[
  \bx = (x,y) = (x,x_{d+1}), \quad x \in \R^d,\ y \in \R^+.
\]

%--------------------------------------------------------------------------------------
% \subsection{The Stinga-Torrea extension}
\subsection{Extension Property}
\label{sub:CaffarelliSilvestre}
%--------------------------------------------------------------------------------------

The groundbreaking extension \eqref{alpha_harm_intro}
of Caffarelli and Silvestre \cite{CS:07}, valid for any power $s \in (0,1)$,
is formulated in $\R^d$. Cabr\'e and Tan \cite{CT:10} and Stinga and Torrea \cite{ST:10} realized that a similar
extension holds for the spectral Laplacian over $\Omega$ bounded; see also \cite{CDDS:11,BCdPS:12}. This problem reads
\begin{equation}
\label{alpha_harm}
  \DIV\left( y^\alpha \nabla \ue \right) = 0  \text{ in } \C, \quad
  \ue = 0  \text{ on } \partial_L \C, \quad
  \partial_{\nu^\alpha} \ue = d_s f  \text{ on } \Omega \times \{0\},
\end{equation}
where $\alpha = 1-2s$ and the so-called conormal exterior derivative of $\ue$ at $\Omega \times \{ 0 \}$ is
\begin{equation}
\label{def:lf}
\partial_{\nu^\alpha} \ue = -\lim_{y \rightarrow 0^+} y^\alpha \partial_y \ue(\cdot,y).
\end{equation}
The limit in \eqref{def:lf} must be understood in the distributional
sense \cite{CT:10,CS:07,CDDS:11,ST:10}. With this construction at hand, the fractional Laplacian and the Dirichlet-to-Neumann operator of problem \eqref{alpha_harm} are related by
\begin{equation*}
%\label{eq:identity} %(JPB- we didn't make reference to this eqn)
  d_s \Laps u = \partial_{\nu^\alpha} \ue \quad \text{in } \Omega.
\end{equation*}

The operator in \eqref{alpha_harm} is in divergence form and thus
amenable to variational techniques. However, it is \emph{nonuniformly}
elliptic because the weight $y^{\alpha}$ either blows up for $-1<\alpha<0$
or degenerates for $0<\alpha<1$ as $y \downarrow0$; the exceptional case $\alpha=0$
corresponds to the regular harmonic extension for $s=\frac12$ \cite{CT:10}. This entails
dealing with weighted Lebesgue and Sobolev spaces with the 
weight $|y|^{\alpha}$ for $\alpha \in (-1,1)$ \cite{BCdPS:12,CT:10,CS:07, CDDS:11}. 
Such a weight belongs to the Muckenhoupt class $A_2(\R^{d+1})$, which
is the collection of weights $\omega$ so that
\cite{Javier,FKS:82,GU,Muckenhoupt,Turesson}
\begin{equation*}
%  \label{A_pclass} %JPB
  C_{2,\omega} = \sup_{B} \left( \fint_{B} \omega \diff \bx \right)
            \left( \fint_{B} \omega^{-1} \diff \bx \right) < \infty,
\end{equation*}
where the supremum is taken over all balls $B$ in $\R^{d+1}$ and
$\fint_B$ stands for the mean value over $B$. The Muckenhoupt
characteristic $C_{2,\omega}$ appears in all estimates involving $\omega$.

If $D \subset \R^d\times\R^+$, we define $L^2(y^\alpha,D)$ as the Lebesgue space 
for the measure $y^{\alpha}\diff \bx$. We also define the weighted Sobolev space
\[
H^1(y^{\alpha},D) =
  \left\{ w \in L^2(y^{\alpha},D): | \nabla w | \in L^2(y^{\alpha},D) \right\},
\]
where $\nabla w$ is the distributional gradient of $w$. 
We equip $H^1(y^{\alpha},D)$ with the norm
\begin{equation}
\label{wH1norm}
\| w \|_{H^1(y^{\alpha},D)} =
\left(  \| w \|^2_{L^2(y^{\alpha},D)} + \| \nabla w \|^2_{L^2(y^{\alpha},D)} \right)^{\frac{1}{2}}.
\end{equation}
The space $H^1(y^{\alpha},D)$ is Hilbert with the norm
\eqref{wH1norm} and $C^{\infty}(D) \cap H^1(y^{\alpha},D)$ is dense
in $H^1(y^{\alpha},D)$ because $|y|^{\alpha} \in A_2(\R^{d+1})$
(cf.~\cite[Theorem~1]{GU}, \cite{KO84} and \cite[Proposition 2.1.2, Corollary 2.1.6]{Turesson}).

To analyze problem \eqref{alpha_harm} we define the weighted Sobolev space
\begin{equation*}
%  \label{HL10} %JPB
  \HL(y^{\alpha},\C) = \left\{ w \in H^1(y^\alpha,\C): w = 0 \textrm{ on } \partial_L \C\right\}.
\end{equation*}
The following \emph{weighted Poincar\'e inequality} holds \cite[inequality (2.21)]{NOS}
\begin{equation*}
%\label{Poincare_ineq} %JPB
\| w \|_{L^2(y^{\alpha},\C)} \lesssim \| \nabla w \|_{L^2(y^{\alpha},\C)}
\quad \forall w \in \HL(y^{\alpha},\C).
\end{equation*}
Consequently, the seminorm on $\HL(y^{\alpha},\C)$ is equivalent to \eqref{wH1norm}. For $w \in H^1(y^{\alpha},\C)$, $\tr w$ denotes its trace onto $\Omega \times \{ 0 \}$ which satisfies
\cite[Proposition 2.5]{NOS}
\begin{equation*}
%\label{Trace_estimate} %JPB
\tr \HL(y^\alpha,\C) = \Hs,
\qquad
  \|\tr w\|_{\Hs} \leq C_{\tr} \| w \|_{\HLn(y^\alpha,\C)}.
\end{equation*}

The variational formulation of \eqref{alpha_harm} reads: find $\ue \in \HL(y^{\alpha}, \C)$ such that
\begin{equation}\label{alpha_harm_weak}
 \int_{\C} y^\alpha \nabla \ue \cdot \nabla w \diff \bx
  = d_s \langle f, \tr w \rangle \quad \forall w \in \HL(y^{\alpha},\C),
\end{equation}
where, we recall that, $\langle \cdot, \cdot \rangle$ corresponds to
the duality pairing between $\Hsd$ and $\Hs$.
The fundamental result of Caffarelli and Silvestre \cite{CS:07} for
$\Omega=\rn$ and of Cabr\'e and Tan \cite[Proposition 2.2]{CT:10} and Stinga and Torrea \cite[Theorem 1.1]{ST:10} for $\Omega$ bounded reads: given $f \in \Hsd$,
if $u \in \Hs$ solves \eqref{fl=f_bdddom} and
$\ue \in \HL(y^{\alpha},\C)$ solves \eqref{alpha_harm_weak},
 then  
\begin{equation}
\label{eq:identity2}
u = \tr \ue,\quad
d_s \Laps u = \partial_{\nu^\alpha} \ue \qquad \text{in } \Omega,
\end{equation}
where the first equality holds in $\Hs$, whereas the second one in $\Hsd$.

%---------------------------------------------------------------------------------------
\subsection{Regularity}\label{S:reg-spectral}
%---------------------------------------------------------------------------------------

To study the finite element discretization of \eqref{alpha_harm_weak} we must understand the regularity of $\ue$. We begin by recalling that if $u = \sum_{k=1}^{\infty} u_k \varphi_k $ solves \eqref{fl=f_bdddom}, then $u_k = \lambda_k^{-s} f_k$, with $f_k$ being the $k$-th Fourier coefficient of $f$. The unique solution $\ue$ of problem \eqref{alpha_harm_intro} thus admits the representation \cite[formula (2.24)]{NOS}
\begin{equation*} 
%\label{eq:SepVar} %JPB
\ue(x,y) = \sum_{k=1}^\infty u_k \varphi_k(x) \psi_k(y).
\end{equation*}
The functions $\psi_k$ solve the 2-point boundary value problem in $\R^+$
\begin{equation*}
%\label{psik} %JPB
\psi_k'' + \frac{\alpha}{y}\psi_k'= \lambda_k \psi_k,
\quad \text{in } (0,\infty);
\qquad \psi_k(0) = 1, \quad \lim_{y\rightarrow \infty} \psi_k(y) = 0.
\end{equation*}
Thus, if $s = \tfrac{1}{2}$, we have that $\psi_k(y) = \exp(-\sqrt{\lambda_k}y)$ \cite[Lemma 2.10]{CT:10}; and, if $s \in (0,1) \setminus \{ \tfrac{1}{2}\}$, then \cite[Proposition 2.1]{CDDS:11}
\[
 \psi_k(y) = c_s (\sqrt{\lambda_k}y)^s K_s(\sqrt{\lambda_k}y),
\]
where $c_s = 2^{1-s}/\Gamma(s)$ and $K_s$ denotes the modified Bessel
function of the second kind \cite[Chap. 9.6]{Abra}. Using asymptotic
properties of $K_s(\zeta)$ as $\zeta\downarrow0$ \cite[Chapter 9.6]{Abra}, and \cite[Chap. 7.8]{MR1429619}, we obtain for $y\downarrow0$:
\[
\psi_k'(y) \approx y^{-\alpha},
\qquad
\psi_k''(y) \approx y^{-\alpha-1}.
\]
Exploiting these estimates, the following regularity results hold
\cite[Theorem 2.7]{NOS}.

\begin{theorem}[global regularity of the $\alpha$-harmonic extension]
\label{T:regularity_extension}
Let $f \in \Ws$, where $\Ws$ is defined in \eqref{e:tildeHs} for $s\in(0,1)$. Let $\ue \in \HL(y^{\alpha},\C)$ solve \eqref{alpha_harm_intro} with $f$ as data. 
Then, for $s \in (0,1) \setminus \{\sr\}$, we have
\begin{equation}
\label{eq:reginx}
  \| \Delta_{x} \ue \|^2_{L^2(y^{\alpha},\C)} + \| \partial_y \nabla_{x} \ue \|^2_{L^2(y^{\alpha},\C)}
  = d_s \| f \|_{\Ws}^2,
\end{equation}
and
\begin{equation}
\label{eq:reginy}
  \| \partial_{yy} \ue \|_{L^2(y^{\beta},\C)} \lesssim \| f \|_{L^2(\Omega)},
\end{equation}
with $\beta>2\alpha+1$. For the special case $s = \sr$, we obtain
\[
 \| \ue \|_{H^2(\C)} \lesssim \| f \|_{\mathbb{H}^{1/2}(\Omega)}.
\]
\end{theorem}

Comparing \eqref{eq:reginx} and \eqref{eq:reginy}, we realize that the
regularity of $\ue$ is much worse in the extended $d+1$ direction.
Since $\Omega$ is convex, the following elliptic
regularity estimate is valid \cite{Grisvard}:
\begin{equation}
\label{Omega_regular}
\| w \|_{H^2(\Omega)} \lesssim \| \Delta_{x} w \|_{L^2(\Omega)} \quad \forall w \in H^2(\Omega) \cap H^1_0(\Omega).
\end{equation}
This, combined with \eqref{eq:reginx}, yields the following estimate
for the Hessian $D_{x}^2\ue$ in the variable $x\in\Omega$:
\[
\| D_{x}^2 \ue \|_{L^2(y^\alpha, \C)} \lesssim \| f \|_{\mathbb{H}^{1-s}(\Omega)}.
\]
Further regularity estimates in H\"older and Sobolev norms are derived
in \cite{CaffarelliStinga:15}. However, we do not need them for what follows.

%-----------------------------------------------------------------------------------
\subsection{Truncation}\label{S:truncation}
%-----------------------------------------------------------------------------------
Since $\C$ is unbounded, problem \eqref{alpha_harm} cannot be
approximated with standard finite element techniques. 
However, since the solution $\ue$ of problem \eqref{alpha_harm} 
decays exponentially in $y$  \cite[Proposition 3.1]{NOS}, by
truncating $\C$ to $\C_\Y:=\Omega\times(0,\Y)$
and setting a homogeneous Dirichlet condition on $y = \Y$, we only incur in an exponentially 
small error in terms of $\Y$ \cite[Theorem 3.5]{NOS}.
If
\[
  \HL(y^{\alpha},\C_\Y) := \left\{ w \in H^1(y^\alpha,\C_\Y): w = 0 \text{ on }
    \Gamma_D \right\},
\]
where $\Gamma_D = \partial_L \C_{\Y} \cup \Omega \times \{ \Y\}$ is
the Dirichlet boundary, then the aforementioned problem reads:
\begin{equation}
\label{alpha_harmonic_extension_weak_T}
\ve \in \HL(y^{\alpha}, \C_\Y): \quad
\int_{\C_\Y} y^\alpha \nabla \ve \cdot \nabla w \diff \bx
  = d_s \langle f, \tr w \rangle \quad \forall w \in \HL(y^{\alpha},\C_\Y).
\end{equation}
The following exponential decay rate is proved in \cite[Theorem 3.5]{NOS}.

\begin{lemma}[truncation error]\label{L:truncation}
If $\ue$ and $\ve$ denote the solutions of \eqref{alpha_harm_weak} and
\eqref{alpha_harmonic_extension_weak_T}, respectively, then
\begin{equation*}
% \label{le:v-v^TC}
  \| \nabla(\ue - \ve) \|_{L^2(y^{\alpha},\C )} \lesssim e^{-\sqrt{\lambda_1} \Y/4} \| f\|_{\Hsd}
\end{equation*}
is valid, where $\lambda_1$ denotes the first eigenvalue of the Dirichlet
Laplace operator and $\Y$ is the truncation parameter.
\end{lemma}

%---------------------------------------------------------------------------------------
\subsection{FEM: A Priori Error Analysis}\label{S:1-FEM-apriori}
%---------------------------------------------------------------------------------------

The first numerical work that exploits the groundbreaking identity \eqref{eq:identity2},
designs and analyzes a FEM for \eqref{alpha_harm} is \cite{NOS}; 
see also \cite{NOSchina,MR3259962}. We briefly review now the main a priori results of \cite{NOS}. 

We introduce a conforming and shape regular mesh $\T = \{K\}$ of $\Omega$, where $K \subset \R^d$ is an element that is isoparametrically equivalent either to the unit cube $[0,1]^d$ or the unit simplex in $\R^d$. Over this mesh we construct the finite element space
\begin{equation}
\label{eq:defofmathbbU}
  \U(\T) = \left\{ W \in C^0(\bar\Omega): W|_K \in \mathcal{P}_1(K) \ \forall K \in \T, \ W|_{\partial\Omega} = 0 \right\}.
\end{equation}
The set $\mathcal{P}_1(K)$ is either the space $\mathbb{P}_1(K)$ of polynomials of total degree at most $1$, when $K$ is a simplex, or the space of polynomials $\mathbb{Q}_1(K)$ of degree not larger than $1$ in each variable provided $K$ is a $d$-rectangle.

To triangulate the truncated cylinder $\C_\Y$ we consider a partition $\mathcal{I}_\Y$ of the interval $[0,\Y]$ with nodes
$y_m, \, m = 0, \cdots, M$, and construct a mesh $\T_{\Y}$ of $\C_\Y$ as the tensor product of $\T$ and $\mathcal{I}_\Y$.
The set of all triangulations of $\C_\Y$ obtained with this procedure is $\Tr$. 

For $\T_{\Y} \in \Tr$, we define
the finite element space 
\begin{equation}
\label{eq:FESpace}
  \V(\T_\Y) = \left\{
            W \in C^0( \bar{\C}_\Y): W|_T \in \mathcal{P}_1(K) \otimes \mathbb{P}_1(I) \ \forall T \in \T_\Y, \
            W|_{\Gamma_D} = 0
          \right\},
\end{equation}
where $T = K \times I$ and observe that $\U(\T)=\tr \V(\T_{\Y})$. Note that $\#\T_{\Y} = M \, \# \T$, and that $\# \T \approx M^d$ implies $\#\T_\Y \approx M^{{d}+1}$.

The Galerkin approximation of \eqref{alpha_harmonic_extension_weak_T} is the function 
$V \in \V(\T_{\Y})$ such that
\begin{equation}
\label{harmonic_extension_weak}
\int_{\C_\Y} y^{\alpha}\nabla V \cdot \nabla W \diff \bx
= 
d_s \langle f, \textrm{tr}_{\Omega} W \rangle \quad \forall W \in \V(\T_{\Y}).
\end{equation}
Existence and uniqueness of $V$ immediately follows from 
$\V(\T_\Y) \subset \HL(y^{\alpha},\C_\Y)$ and the Lax-Milgram Lemma.
It is trivial also to obtain a best approximation result \emph{\`a la}
C\'ea, namely
\begin{equation}\label{best-approx-spectral}
  \| \nabla(\ve - V ) \|_{L^2(y^\alpha, \C_{\Y})}
  = \inf_{W\in\V(\T_{\Y})} \| \nabla(\ve - W ) \|_{L^2(y^\alpha, \C_{\Y})}.
\end{equation}  
This reduces the numerical analysis of \eqref{harmonic_extension_weak}
to a question in approximation theory, which in turn can be answered by
the study of piecewise polynomial interpolation in Muckenhoupt
weighted Sobolev spaces; see \cite{NOS,NOS2}. Exploiting the Cartesian
structure of the mesh $\T_{\Y}$ we are able to extend the anisotropic
estimates of Dur\'an and Lombardi \cite{DL:05} to our setting
\cite[Theorems 4.6--4.8]{NOS}, \cite{NOS2}. The following error
estimates separate in each direction.

\begin{proposition}[anisotropic interpolation estimates]\label{P:anisotropic-est}
There exists a quasi interpolation operator $\Pi_{\T_\Y} : L^1(\C_\Y)\to \V(\T_\Y)$
that satisfies the following anisotropic error estimates for all
$j=1,\ldots,d+1$ and all $T = K \times I \in\T_\Y$
\begin{align*}
  \| w - \Pi_{\T_\Y} w \|_{L^2(y^\alpha,T)} & \lesssim 
    h_K  \| \nabla_{x} w \|_{L^2(y^\alpha,S_T)} + h_I \| \partial_y w \|_{L^2(y^\alpha,S_T)}, 
  \\
  \| \partial_{x_j}(w - \Pi_{\T_\Y} w) \|_{L^2(y^\alpha,T)} 
  &\lesssim
   h_K  \| \nabla_{x} \partial_{x_j} w\|_{L^2(y^\alpha,S_T)} 
 + h_I \| \partial_y \partial_{x_j} w\|_{L^2(y^\alpha,S_T)},
\end{align*}
where $S_T$ stands for the patch of elements of
$\T_\Y$ that intersect $T$, $h_K = |K|^{1/d}$ and $h_I = |I|$.
\end{proposition}

As a first application of Proposition \ref{P:anisotropic-est} we consider a quasiuniform mesh $\T_\Y$ of size $h$ and set  $w = \ue$.
% Proposition \ref{P:anisotropic-est} to $v=\ue$ over a quasi-uniform
% mesh $\T_{\Y}$ of size $h_\T$.
We estimate $\ue-\Pi_{\T_\Y}\ue$ for $y\ge 2h$ as follows:
\[
\int_{2h}^\Y y^\alpha \|\partial_y (\ue-\Pi_{\T_\Y}\ue)\|^2_{L^2(\Omega)} \diff y
\lesssim h^2 \int_h^\Y y^\alpha \left( \|\partial_{yy}\ue\|_{L^2(\Omega)}^2
+ \|\nabla_{x}\partial_y\ue\|_{L^2(\Omega)}^2 \right) \diff y.
\] 
For the first term we resort to \eqref{eq:reginy}, and recall that
$\beta > 2\alpha+1 >\alpha$, to deduce
\begin{align*}
h^2 \int_h^\Y y^\alpha \|\partial_{yy}\ue\|_{L^2(\Omega)}^2 \diff y &\le 
h^2 \sup_{h \le y\le \Y} y^{\alpha-\beta} \int_0^\Y y^\beta 
\|\partial_{yy}\ue\|_{L^2(\Omega)}^2 \diff y \\
&\le h^{2+\alpha-\beta} \|f\|_{L^2(\Omega)}^2.
\end{align*}
For the second term we use \eqref{eq:reginx} instead to arrive at
\[
h^2 \int_h^\Y y^\alpha \|\nabla_{x}\partial_y\ue\|_{L^2(\Omega)}^2 \diff y
\le h^2 \int_0^\Y y^\alpha \|\nabla_{x}\partial_y\ue\|_{L^2(\Omega)}^2 \diff y
\lesssim h^2 \|f\|_{\Ws}^2.
\]
Combining the two estimates we derive the interpolation estimate
\[
\int_{2h}^\Y y^\alpha \|\partial_y(\ue-\Pi_{\T_\Y}\ue)\|^2_{L^2(\Omega)} \diff y
\lesssim h^{2+\alpha-\beta} \|f\|_{\Ws}^2,
\]
which is quasi-optimal in terms of regularity because
$2+\alpha-\beta = 2(s-\epsilon)$ and 
$\partial_y\ue \approx y^{-\alpha}$ formally implies $\ue \in H^{s-\epsilon}(y^\alpha,\C)$
for any $\epsilon>0$; however this estimate exhibits a
suboptimal rate. To restore a quasi-optimal
rate, we must compensate the behavior of $\partial_{yy}\ue$ by a
\emph{graded} mesh in the extended direction, which is allowed by
Proposition \ref{P:anisotropic-est}.
Therefore, we construct a mesh $\mathcal{I}_\Y$ with nodes
\begin{equation}
\label{graded_mesh}
  y_m = m^\gamma M^{-\gamma} \Y, \quad m=0,\dots,M,
\end{equation}
where $\gamma > 3/(1-\alpha)=3/(2s)$.
Combining \eqref{best-approx-spectral} with Proposition \ref{P:anisotropic-est}
we obtain estimates in terms of degrees of freedom 
\cite[Theorem 5.4 and Corollary 7.11]{NOS}.

\begin{theorem}[a priori error estimate]
\label{TH:fl_error_estimates}
Let $\T_\Y=\T \times \mathcal{I}_\Y \in \Tr$ with $\mathcal{I}_\Y$
satisfying \eqref{graded_mesh}, and let $\V(\T_\Y)$ be defined by \eqref{eq:FESpace}. 
If $V \in \V(\T_\Y)$ solves \eqref{harmonic_extension_weak}, 
we have
\begin{equation*}
\label{optimal_rate}
  \| \nabla(\ue - V) \|_{L^2(y^\alpha,\C)} \lesssim
|\log(\# \T_{\Y})|^s(\# \T_{\Y})^{-1/(d+1)} \|f \|_{\mathbb{H}^{1-s}(\Omega)},
\end{equation*}
where $\Y \approx \log(\# \T_{\Y})$. 
Alternatively, if $u$ denotes the solution of \eqref{fl=f_bdddom}, then
\begin{equation*}
\| u - U \|_{\Hs} \lesssim
|\log(\# \T_{\Y})|^s(\# \T_{\Y})^{-1/(d+1)} \|f \|_{\mathbb{H}^{1-s}(\Omega)},
\end{equation*}
where $U = \tr V$.
\end{theorem}

\begin{remark}[domain and data regularity]\label{R:dom_and_data2}
The estimates of Theorem~\ref{TH:fl_error_estimates} hold only if $f \in \mathbb{H}^{1-s}(\Omega)$ 
and the domain $\Omega$ is such that \eqref{Omega_regular} is valid.
\end{remark}

\begin{remark}[quasi-uniform meshes]\label{R:quasiuniform}
Let $V \in \V(\T_\Y)$ solve \eqref{harmonic_extension_weak} over a
quasi-uniform mesh $\T_{\Y}$ of $\C_\Y$ of size $h$. Combining
\eqref{best-approx-spectral} with the preceding discussion and
accounting for the missing domain $\Omega\times(0,2h)$, which yields a better
estimate \cite{NOS}, we obtain for
$\Y \approx |\log h |$ and all $\varepsilon >0$
\begin{equation*}
%  \label{sub_optimal} %JPB
   \| \nabla(\ue - V ) \|_{L^2(y^\alpha, \C_{\Y})}  \lesssim { h^{s-\varepsilon} }\| f \|_{\Ws},
\end{equation*}
where the hidden constant blows up if $\varepsilon\downarrow 0$.
\end{remark}

\begin{remark}[complexity]\label{R:complexity}
Except for a logaritmic factor, the error estimates of Theorem
\ref{TH:fl_error_estimates} decay with a rate $(\#
\T_{\Y})^{-1/(d+1)}$, where $d$ is the dimension of $\Omega$. The FEM
\eqref{harmonic_extension_weak} is thus sub-optimal as a method to
compute in $\Omega$. This can be improved to an error decay $(\#\T_{\Y})^{-1/d}$
with geometric grading; see Section \ref{sub:ext_NOS}.
\end{remark}

\begin{remark}[case $s = \srn$]
If $s=\srn$, we obtain the optimal estimate
\[
\| \nabla(\ue - V ) \|_{L^2(\C_{\Y})} 
     \lesssim h \| f \|_{\mathbb{H}^{1/2}(\Omega)}.
\]
\end{remark}

%-----------------------------------------------------------------------------------
\subsection{Numerical Experiments}\label{S:1-FEM-num}
%-----------------------------------------------------------------------------------

We present two numerical examples for $d=2$ computed within the
\texttt{deal.II} library \cite{dealii,dealii2} using graded meshes. Integrals are evaluated with Gaussian quadratures of sufficiently high order and linear systems are solved using CG with ILU preconditioner and the exit criterion being that the $\ell^2$-norm of the residual is less than $10^{-12}$.

%------------------------------------------------------------------------------------
\subsubsection{Square Domain}
\label{subsub:square}
%------------------------------------------------------------------------------------

Let $\Omega = (0,1)^2$. Then
\[
  \varphi_{m,n}(x_1,x_2) = \sin(m \pi x_1)\sin(n \pi x_2),
  \quad
  \lambda_{m,n} = \pi^2 \left( m^2 + n^2 \right),
  \qquad m,n \in \mathbb{N}.
\]
If $f(x_1,x_2) = ( 2\pi^2)^{s} \sin(\pi x_1)\sin(\pi x_2)$, then
$
  u(x_1,x_2) = \sin(\pi x_1)\sin(\pi x_2),
$
by \eqref{def:second_frac} and
$
  \ue(x_1,x_2,y) = 2^{1-s/2} \pi^s \Gamma(s)^{-1} \sin(\pi x_1)\sin(\pi x_2) y^{s}K_s(\sqrt{2}\pi y)
$
\cite[(2.24)]{NOS}.

We construct a sequence of  meshes $\{\T_{\Y_k} \}_{k\geq1}$, where
$\T$ is obtained by uniform refinement and $\calI_{\Y_k}$ is
given by \eqref{graded_mesh} with parameter $\gamma >
3/(1-\alpha)$. On the basis of Theorem~\ref{TH:fl_error_estimates},
the truncation parameter $\Y_k$ is chosen to be
\[
  \Y_k \geq \frac{2}{\sqrt{\lambda_1}} ( \log C - \log (\# \T_{\Y_{k-1}})^{-1/3} ).
\] 
With this type of meshes,
\[
  \| u - \tr V_k \|_{\Hs} \lesssim  \| \ue -V_k \|_{\HLn(y^\alpha,\C)} \lesssim  {|\log( \# \T_{\Y_k} )|^s} \cdot (\# \T_{\Y_k})^{-1/3},
\]
which is near-optimal in $\ue$ but suboptimal in $u$, since we should expect (see \cite{BrennerScott})
\[
  \| u - \tr V_k  \|_{\Hs} \lesssim h^{2-s} \lesssim (\# \T_{\Y_k} )^{-(2-s)/3}.
\] 

\begin{figure}[h!]
\label{fig:02}
\begin{center}
  \includegraphics[scale=0.3]{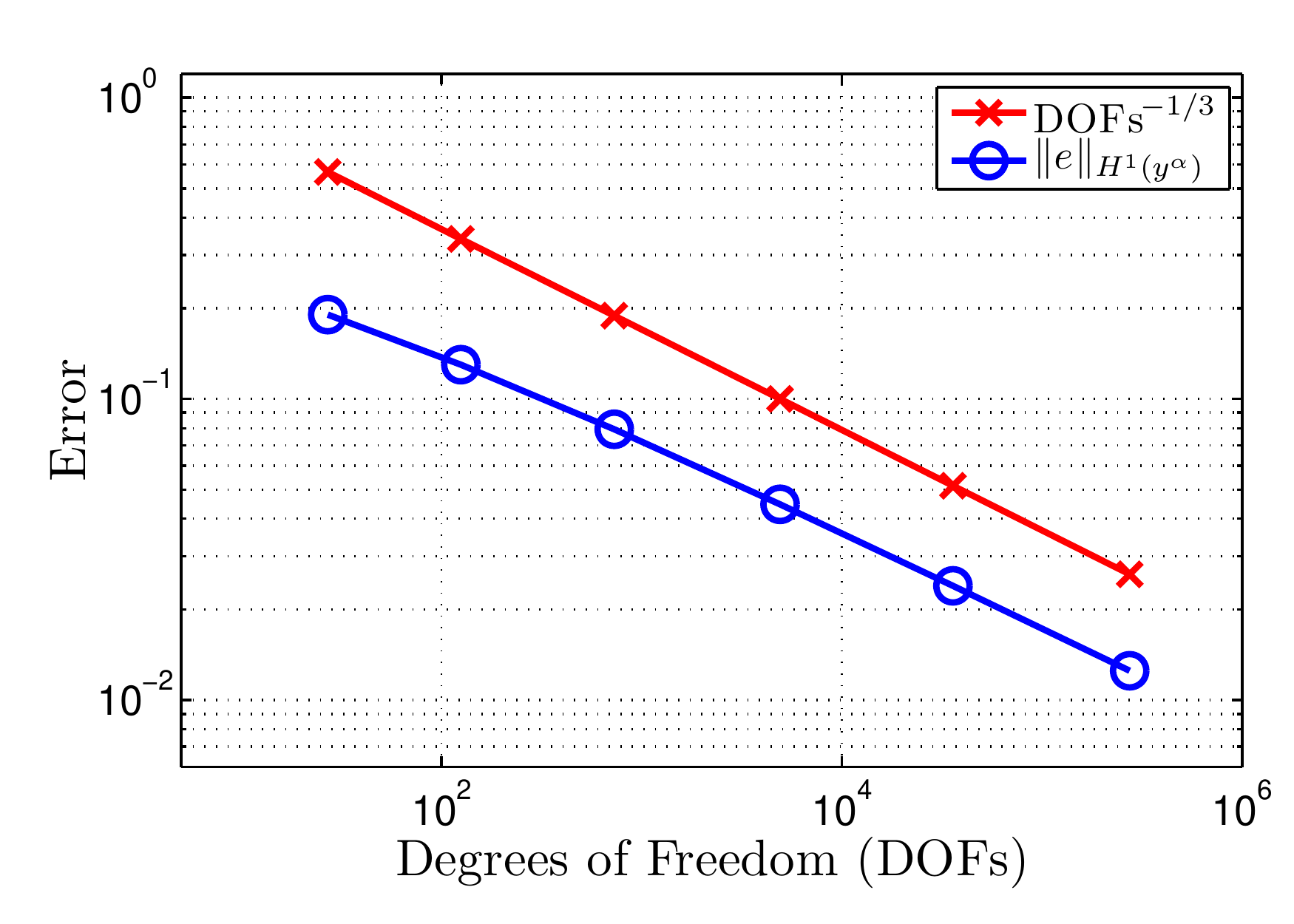}
  \includegraphics[scale=0.3]{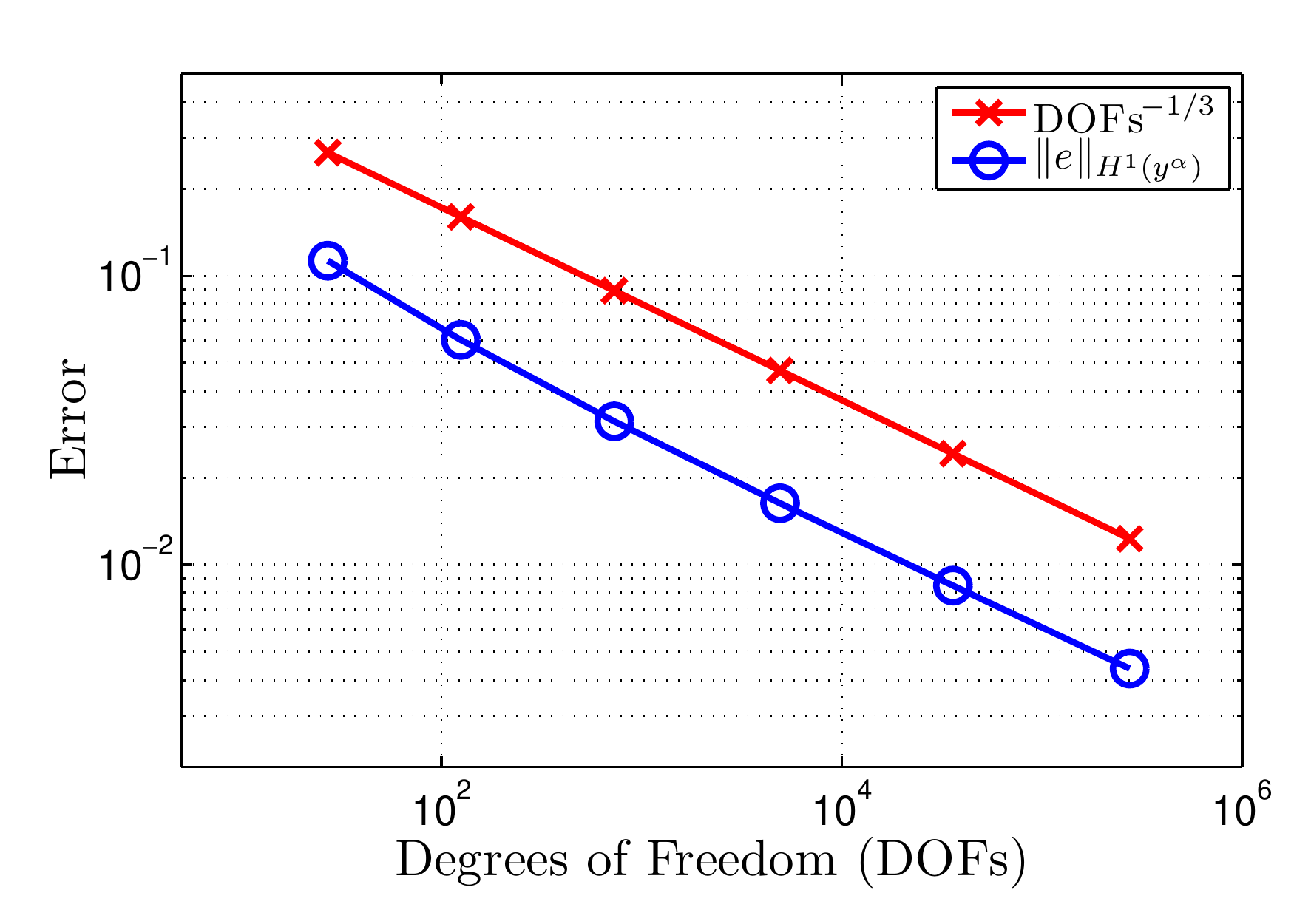}
\end{center}
\vskip-0.3cm
\caption{Computational rate of convergence for a square and graded meshes. The left panel shows the rate for $s=0.2$ and the right one 
for $s=0.8$. The experimental rate of convergence is $(\# \T_{\Y})^{-1/3}$, in agreement with Theorem~\ref{TH:fl_error_estimates}.}
\end{figure}

Figure~\ref{fig:02} shows the rates of convergence for $s=0.2$ and $s=0.8$ respectively.
In both cases, we obtain the rate given by Theorem~\ref{TH:fl_error_estimates}.

%------------------------------------------------------------------------------------
\subsubsection{Circular Domain}
%------------------------------------------------------------------------------------

If $\Omega = \{ x \in \R^2 : |x|<1 \}$, then
\begin{equation*}
%\label{circularphi} %JPB
 \varphi_{m,n}(r,\theta) = 
 J_{m}(j_{m,n}r)
\left( A_{m,n} \cos (m \theta) + B_{m,n} \sin(m \theta)\right),
\
\lambda_{m,n}=j_{m,n}^2,
\end{equation*}
where $J_m$ is the $m$-th Bessel function of the first kind;
$j_{m,n}$ is the $n$-th zero of $J_m$ and
$A_{m,n}$, $B_{m,n}$ are normalization constants to ensure
$\| \varphi_{m,n} \|_{L^2(\Omega)}=1$.

% From \cite[Chapter 9]{Abra}, we have that $j_{1,1} \approx 3.8317$. Thus, 
If $f = ( \lambda_{1,1})^{s} \varphi_{1,1}$, then \eqref{def:second_frac} and \cite[(2.24)]{NOS} show that $u = \varphi_{1,1}$ and
\[
  \ue(r,\theta,y) = 2^{1-s}\Gamma(s)^{-1}(\lambda_{1,1})^{s/2}
         \varphi_{1,1}(r,\theta) y^{s}K_s(\sqrt{2}\pi y).
\]

\begin{figure}
\begin{center}
  \includegraphics[scale=0.3]{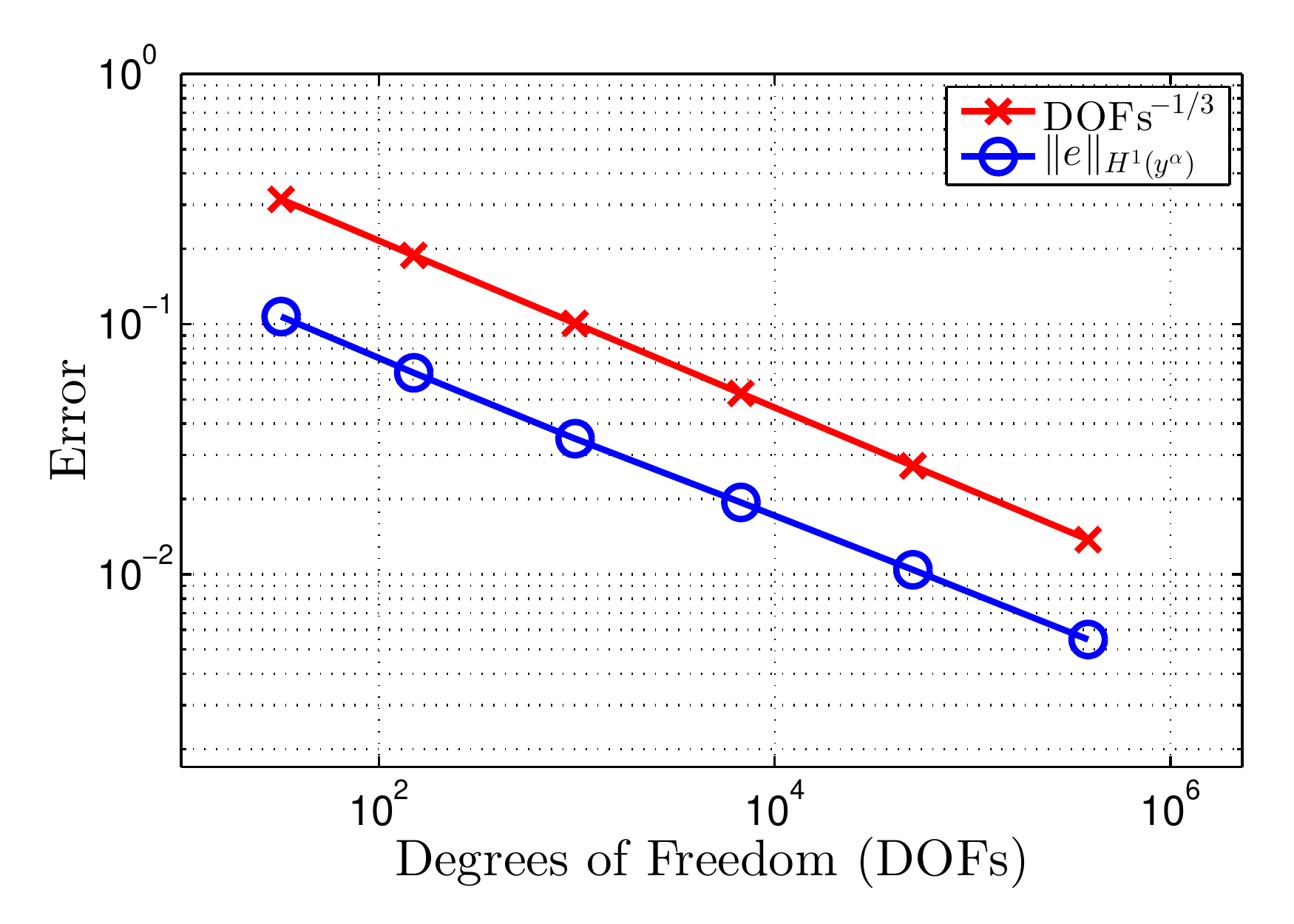}
  \includegraphics[scale=0.3]{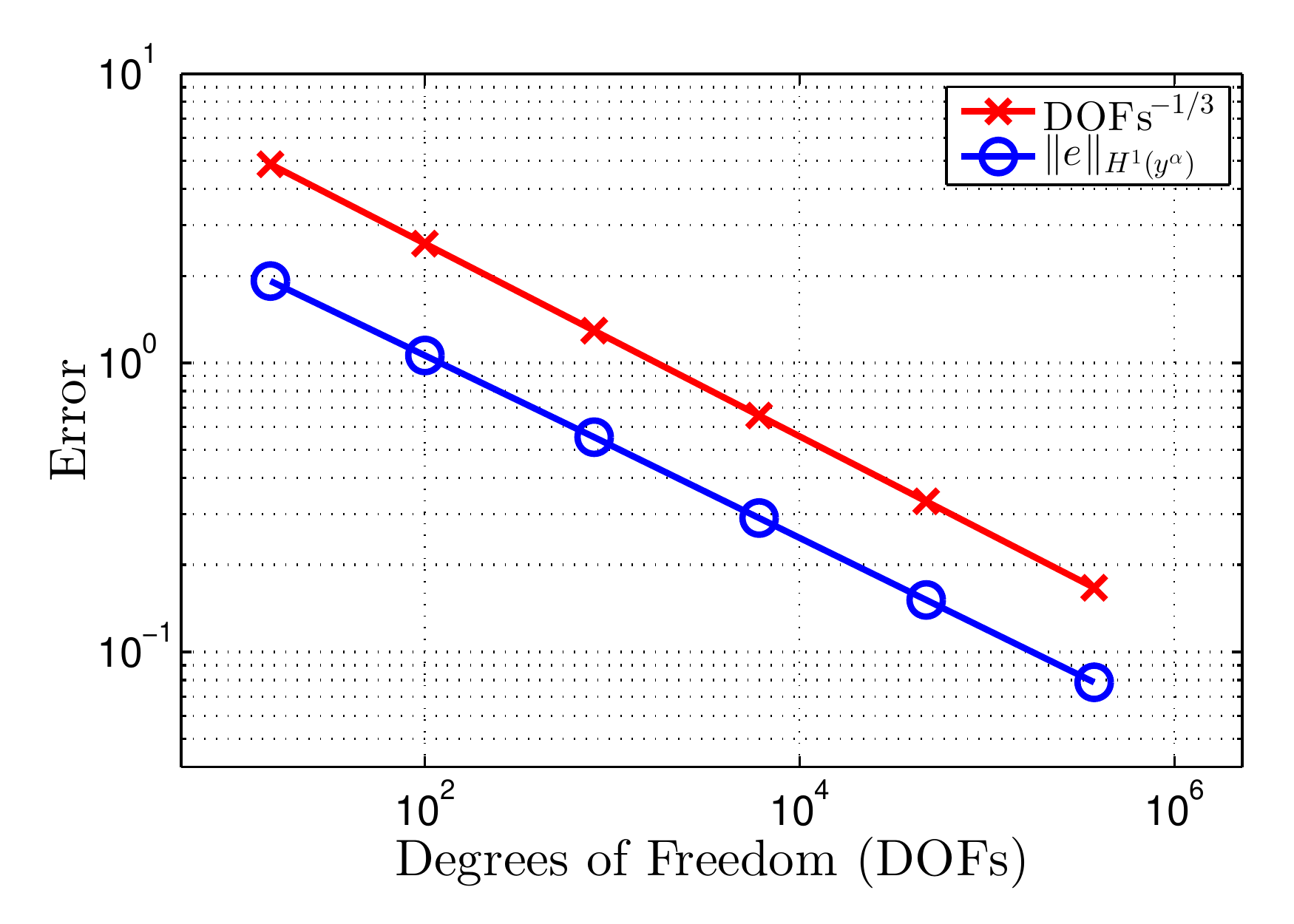}
\end{center}
\vskip-0.3cm
\caption{Computational rate of convergence for a circle with graded meshes. The left panel shows the rate for $s=0.3$ and the right one 
for $s=0.7$. The experimental rate of convergence is $(\#\T_{\Y} )^{-1/3}$,
in agreement with Theorem~\ref{TH:fl_error_estimates}.}
\label{fig:03}
\end{figure}

We construct a sequence of  meshes $\{\T_{\Y_k} \}_{k\geq1}$ as in \S\ref{subsub:square}. With these meshes
\begin{equation*}
%\label{numerical_experiment_2_ve} %JPB
  \| \ue -V_k \|_{\HLn(y^\alpha,\C)} \lesssim 
  |\log(\# \T_{\Y_k})|^s (\# \T_{\Y_k})^{-1/3},
\end{equation*}
which is near-optimal. Figure~\ref{fig:03} shows the errors of $\|\ue - V_k \|_{H^1(y^{\alpha},\C_{\Y_{k}})}$
for $s = 0.3$ and $s = 0.7$. The results, again, are in agreement with
Theorem~\ref{TH:fl_error_estimates}.

%---------------------------------------------------------------------------------------
\subsection{FEM: A Posteriori Error Analysis}\label{S:1-FEM-aposteriori}
%---------------------------------------------------------------------------------------

A posteriori error estimation and adaptive finite element methods (AFEMs) have been the subject of
intense research since the late 1970's because they yield optimal
performance in situations where classical FEM cannot.
The a priori theory for \eqref{harmonic_extension_weak} requires $f \in \Ws$ and
$\Omega$ convex for \eqref{Omega_regular} to be valid; see Remark \ref{R:dom_and_data2}.
If either of these does not hold, then $\ue$ may have singularities in
the $x$-variables and Theorem~\ref{TH:fl_error_estimates} may not
apply: a quasi-uniform refinement of $\Omega$ would not result in an
efficient solution technique. An adaptive loop driven by an a
posteriori error estimator is essential to recover optimal
rates of convergence. In what follows we explore this.

We first observe that we cannot rely on residual error estimators.
In fact, they hinge on the strong form of the local residual to measure the error. Let $T \in \T_{\Y}$ and $\nu$ be the unit outer normal to $T$. Integration by parts yields
\[
 \int_{T} y^{\alpha} \nabla V \cdot \nabla W \diff \bx = \int_{\partial T}  y^{\alpha} W \nabla V \cdot \nu \diff \sigma
- \int_{T} \DIV(y^{\alpha} \nabla V) W \diff \bx.
\]
Since $\alpha \in (-1,1)$ the boundary integral is meaningless for $y = 0$. Even the very first step in the derivation of a residual a posteriori error estimator fails! There is nothing left to do but to consider a different type of estimator.

%----------------------------------------------------------------------------------
\subsubsection{Local Problems over Cylindrical Stars}\label{S:cylstars}
%----------------------------------------------------------------------------------

Inspired by \cite{BabuskaMiller,MNS02}, we deal with the anisotropy of the mesh in the extended variable $y$ and the coefficient $y^{\alpha}$ by considering local problems on cylindrical stars. The solutions of these local problems allow us to define an anisotropic a posteriori error estimator which, under certain assumptions, is equivalent to the error up to data oscillation terms.

Given a node $\mathbf{v}$ on the mesh $\T_{\Y}$, we exploit the tensor product structure of $\T_{\Y}$, and  we write $\mathbf{v} = (\vero,\texttt{w})$ where $\vero$ and $\texttt{w}$ are nodes on the meshes 
$\T$ and $\mathcal{I}_{\Y}$ respectively. For $K \in \T$, we denote by $\N(K)$ and $\Nin(K)$ the set of nodes and interior nodes of $K$, respectively. We set 
\[
\N(\T) = \bigcup_{K \in \T} \N(K), \qquad \Nin(\T) = \bigcup_{K \in \T} \Nin(K).
\]
The \emph{star} around $\vero$ is
$
  S_{\vero} = \bigcup_{K \ni \vero} K \subset \Omega,
$
and the \emph{cylindrical star} around $\vero$ is
\[
  \C_{\vero} := \bigcup\left\{ T \in \T_\Y : T = K \times I,\ K \ni \vero  \right\}= S_{\vero} \times (0,\Y) \subset \C_{\Y}.
\]

For each node $\vero \in \N(\T)$ we define the local space
\begin{align*}
  \mathbb{W}(\C_{\vero}) &= 
    \left\{
     w \in H^1(y^\alpha,\C_{\vero}): w = 0 ~\textrm{on} ~ \partial
     \C_{\vero} \backslash \Omega\times \{ 0 \}
    \right\},
\end{align*}
and the (ideal) estimator $\eta_{\vero} \in \mathbb{W}(\C_{\vero})$ to be the solution of
\begin{equation*}
  \int_{\C_{\vero}} y^{\alpha} \nabla\eta_{\vero} \cdot \nabla w \diff \bx=
  d_s \langle f, \tr w \rangle
   -  \int_{\C_{\vero}} y^{\alpha} \nabla V \cdot \nabla w \diff \bx
   \quad \forall w \in \mathbb{W}(\C_{\vero}).
\end{equation*}
We finally define the local indicators $\E_{\vero}$ and global error estimators
$\E_{\T_{\Y}}$ as follows:
\begin{equation}
\label{ideal_global_estimator}
  \E_{\vero} = \| \nabla \eta_{\vero} \|_{L^2(y^\alpha,\C_{\vero})}, \qquad
  \E_{\T_{\Y}}^2 = \sum_{\vero \in \N(\T) } \E_{\vero}^2.
\end{equation}
We have the following key properties \cite[Proposition 5.14]{CNOS2}.

\begin{proposition}[a posteriori error estimates]\label{P:dis_lower_bound}
Let $\ve \in \HL(y^\alpha,\C_\Y)$ and $V \in \V(\T_\Y)$ solve 
\eqref{alpha_harmonic_extension_weak_T} and \eqref{harmonic_extension_weak} respectively. 
Then, the estimator defined in \eqref{ideal_global_estimator}
satisfies the global bound
\[
\| \nabla(\ve-V)\|_{L^2(y^{\alpha},\C_\Y)} \lesssim \E_{\T_{\Y}},
\]
and the local bound with constant $1$ for any $\textup{\vero} \in \N(\T_{\Omega})$ 
\begin{equation*}
\E_{\textup{\vero}} \leq\| \nabla(\ve-V)\|_{L^2(y^{\alpha},\C_{\textup{\vero}})}.
\end{equation*}
\end{proposition}

These estimates provide the best scenario for a posteriori error
analysis but they are not practical because the local space
$\mathbb{W}(\C_{\vero})$ is infinite dimensional. We further discretize
$\mathbb{W}(\C_{\vero})$ with continuous piecewise polynomials of degree $>1$ as follows.
If $K$ is a quadrilateral, we use polynomials of degree $\le 2$
in each variable. If $K$ is a simplex, we employ polynomials of
total degree $\le 2$ augmented by a local cubic bubble function. We
tensorize these spaces with continuous piecewise quadratics in the
extended variable. We next construct discrete subspaces of the local
space $\mathbb{W}(\C_{\vero})$ and corresponding discrete estimators
instead of \eqref{ideal_global_estimator}. Under suitable assumptions,
Proposition \ref{P:dis_lower_bound} extends to this case
\cite[Section 5.4]{CNOS2}.

%--------------------------------------------------------------------------------------
\subsubsection{Numerical Experiment}
%--------------------------------------------------------------------------------------

We illustrate the performance of a practical version of the a posteriori error estimator
\eqref{ideal_global_estimator}. We use an adaptive loop
\begin{equation*}
 \textsf{\textup{SOLVE}} \rightarrow \textsf{\textup{ESTIMATE}} \rightarrow \textsf{\textup{MARK}} \rightarrow \textsf{\textup{REFINE}}
%\label{afem} %JPB
\end{equation*}
with D\"orfler marking.
We generate a new mesh $\T'$ by bisecting all the elements 
$K \in \T$ contained in the marked set $\mathscr{M}$ based on newest-vertex
bisection method; \cite{NSV:09,NV}. 
We choose the truncation parameter as $\Y = 1 + \tfrac{1}{3}\log(\# \T')$ \cite[Remark 5.5]{NOS}.
We set $M \approx (\# \T')^{1/d}$ and construct $\mathcal{I}_\Y'$ by the rule \eqref{graded_mesh}. The new mesh
$
  \T_\Y' = \textsf{\textup{REFINE}}(\mathscr{M})
$
is obtained as the tensor product of $\T'$ and $\mathcal{I}_\Y'$.

\begin{figure}[h!]
\centering
\includegraphics[width=0.4\textwidth]{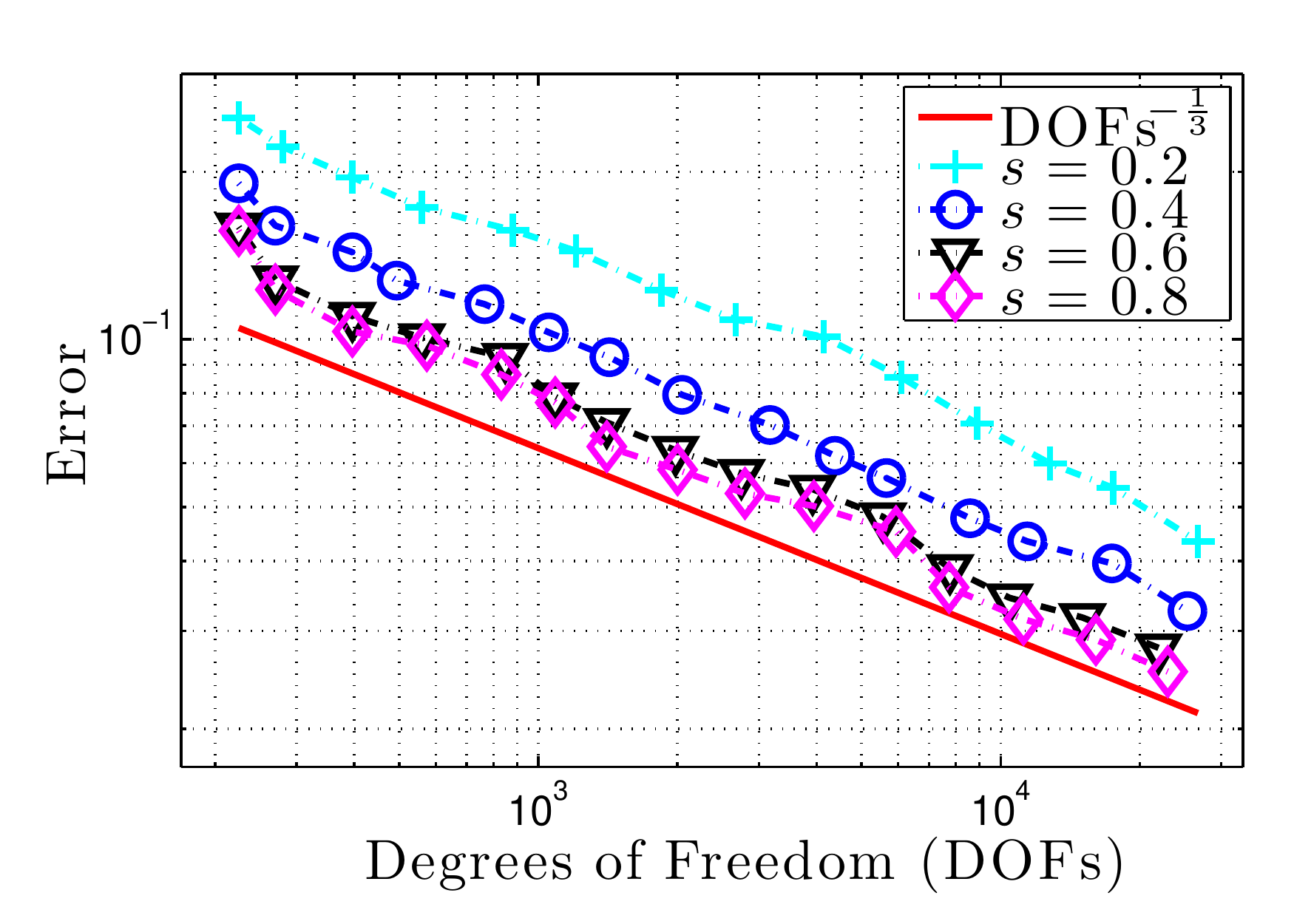}
\hfil
\includegraphics[width=0.4\textwidth]{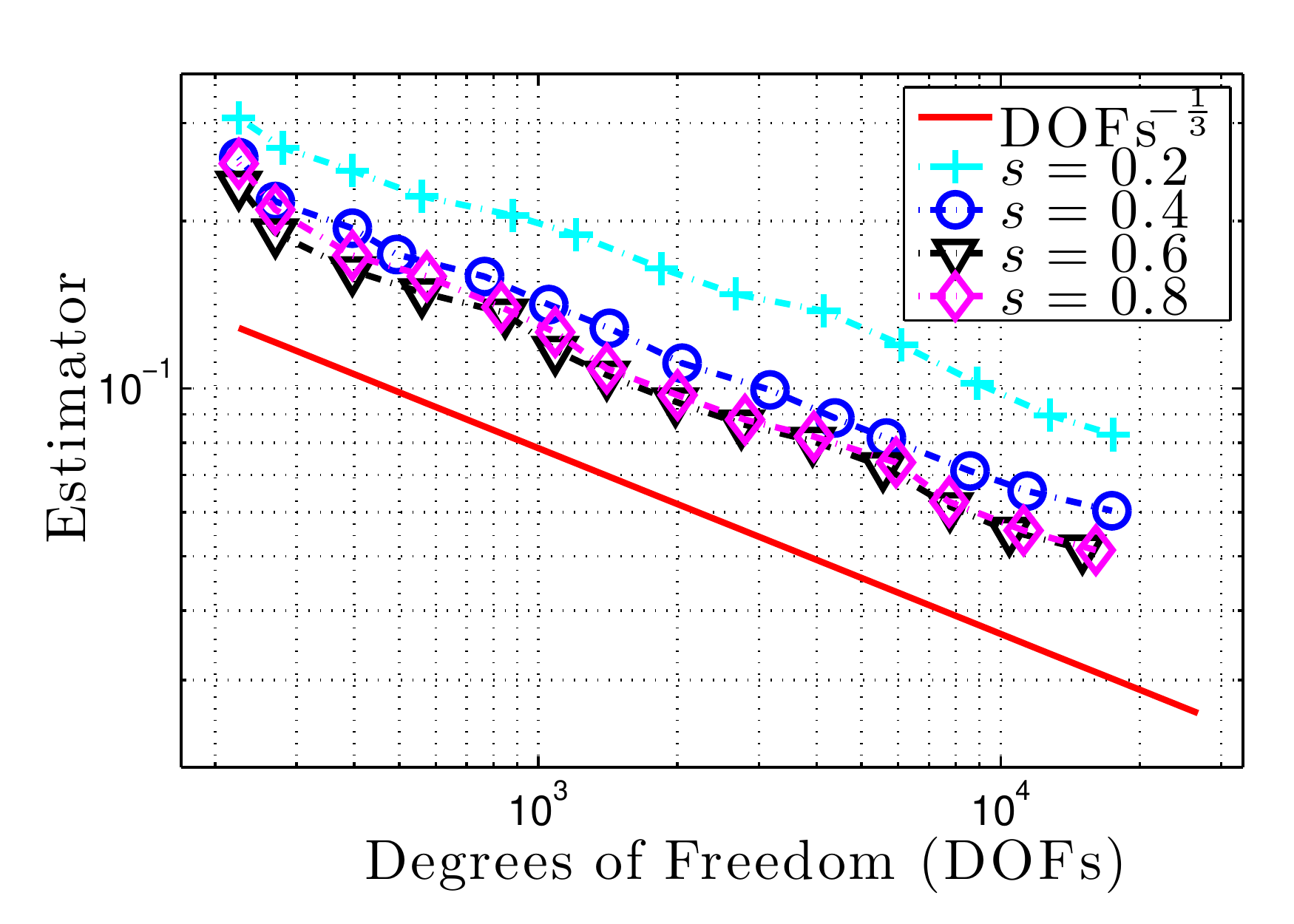}
\caption{Computational rate of convergence of AFEM for $d=2$ and
  $s=0.2$, $0.4$, $0.6$ and $s=0.8$. The left panel shows the 
  error vs.~the number of degrees of freedom, the right one the total
  error indicator. We recover the optimal rate $\#
  (\T_{\Y})^{-1/3}$. For $s<\tfrac12$, the right hand side $f\equiv1\notin
  \Ws$ and a quasiuniform mesh in $\Omega$ does not deliver an optimal
  rate of convergence \cite[Section 6.3]{NOS}.}
\label{F:incompatible}
\end{figure}

We consider the data $\Omega = (0,1)^2$ and $f\equiv1$, which is
incompatible for $s<1/2$ because $f$ does not have a vanishing trace
whence $f\notin \Ws$; see also \cite[Section 6.3]{NOS}. Therefore, we cannot
expect an optimal rate for quasi-uniform meshes $\T$ in $\Omega$
according to Theorem~\ref{TH:fl_error_estimates}. Figure
\ref{F:incompatible} shows that adaptive mesh refinement guided by AFEM
restores an optimal decay rate for $s < 1/2$.

%-------------------------------------------------------------------------------------
\subsection{Extensions and Applications}\label{sub:ext_NOS}
%-------------------------------------------------------------------------------------

We conclude the discussion of this approach by mentioning several extension and applications:

\begin{enumerate}[$\bullet$]
  \item \emph{Efficient solvers}: Finding the solution to \eqref{alpha_harmonic_extension_weak_T} entails solving a large linear system with a sparse matrix. In \cite{CNOS} the construction of efficient multilevel techniques for the solution of this problem was addressed. It was shown that multilevel techniques with line smoothers over vertical lines in the extended direction perform almost optimally; \ie the contraction factor depends linearly on the number of levels, and thus logarithmically on the problem size.
  
  \item \emph{Time dependent problems}: In \cite{NOS3} the time
dependent problem
\begin{equation} \label{eq:CaPUTOivbp}
\partial_t^\gamma u + \Laps u = f\quad \Omega\times(0,T),
\qquad
u(\cdot,0)=u_0
\end{equation}
was examined. Here $\partial_t^\gamma$ denotes the so-called Caputo
derivative of order $\gamma \in (0,1)$, which is defined as follows \cite{MR1347689}:
\begin{equation}
\label{eq:caputo} 
\partial_t^\gamma w (x,t) :=\frac{1}{\Gamma(1-\gamma)} \int_0^t
\frac{1}{(t-r)^\gamma} \partial_r w (x,r) \diff r
\quad t>0, \; x\in\Omega.
\end{equation}
It turns out that the solution of this problem is always singular as
$t\downarrow0$ provided the initial condition $u_0\ne0$. In fact,
taking $f=0$ and representing $u$ in terms of the Mittag-Leffler function \cite{Gorenflo} reveals that
\begin{equation*}
u(x,t) = \left(1 - \frac{t^\gamma}{\Gamma(1+\gamma)} (-\Delta)^s +
  \mathcal{O}(t^{2\gamma})\right) u_0(x)
  \quad x\in\Omega.
\end{equation*}
These heuristics led to the following regularity results shown in \cite{NOS3}
\begin{equation}
\label{eq:regstate}
 \partial_t u \in L\log L(0,T;\Hsd), \qquad \partial_{tt}^2 u \in L^2(t^\sigma;0,T;\Hsd), \ \sigma > 3-2\gamma,
\end{equation}
where $L\log L(0,T)$ denotes the Orlicz space of functions $w$ such that $|w|\log|w| \in L^1(0,T)$; see \cite{MR0126722}.
      Using the extension property, problem \eqref{eq:CaPUTOivbp} reduces to a
      quasi-stationary elliptic equation with a dynamic boundary
      condition. Rates of convergence for fully discrete schemes were
      derived in \cite{NOS3}, which are consistent with the regularity
      \eqref{eq:regstate};
      this issue has been largely ignored in the literature.
      
      The extension of these results to a space-time fractional wave equation, \ie $\gamma \in (1,2]$ is currently under investigation \cite{OS:17}.

  \item \emph{Nonlinear problems}: Elliptic and parabolic obstacle
    problems with the spectral fractional Laplacian were considered in
    \cite{MR3393323} and \cite{MR3542012}, respectively. Rates of
    convergence for their FEM approximation were derived, and this
    required a careful combination of Sobolev regularity, as in
    Theorem~\ref{T:regularity_extension}, and H\"older regularity of
    the solution \cite{MR2367025,MR3100955}. In addition, a positivity preserving interpolant that is stable in anisotropic meshes and weighted norms had to be constructed.
  
  \item \emph{PDE constrained optimization and optimal control}:
    Optimal control problems where the state equation is given by
    either a stationary or parabolic equation with a spectral
    fractional Laplacian were studied in \cite{MR3429730} and
    \cite{MR3504977}, respectively. Existence and uniqueness of
    optimal pairs was obtained, as well as their regularity. In both
    cases fully discrete schemes were designed and their convergence
    shown. The results of \cite{MR3429730} were later improved upon
    considering piecewise linear approximation of the optimal control
    variable \cite{Otarola_controlp1} and adaptive algorithms
    \cite{AO2}. Sparse optimal control for the spectral fractional
    Laplacian was studied in \cite{OS2}. Finally, reference
    \cite{AOS3} presents the design and analysis of an approximation
    scheme for an optimal control problem where the control variable
    corresponds to the order of the fractional operator
    \cite{SprekelsV}.

   \item \emph{Near optimal complexity}: According to Remark~\ref{R:complexity},
    the FEM on anisotropic meshes with radical grading
    \eqref{graded_mesh} is suboptimal. This deficiency can be cured by
    exponential grading in the extended variable $y$ and suitable
    exploitation of analyticity properties of $\ue(\cdot,y)$ in
    $y$. This leads to a hybrid FEM that combines the $hp$-version in
    $y$ with the $h$-version in $\Omega$ and exhibits an error decay
    $(\# \T_{\Y})^{-1/d}$ provided $f\in\Ws$
    \cite{BMNOSS:17,MFSV:17}. Dealing with incompatible $f\notin\Ws$
    is open although this is responsible for boundary singularities
    governed by \eqref{boundary-CS}.

\end{enumerate}

% !TEX root = three-methods.tex
%%%%%%%%%%%%%%%%%%%%%%%%%%%%%%%%%%%%%%%%%%%%%%%%%%%%%%%%%%%%%%%%%%%%%%%%%%%%%%%%%%%%

%%%%%%%%%%%%%%%%%%%%%%%%%%%%%%%%%%%%%%%%%%%%%%%%%%%%%%%%%%%%%%%%%%%%%%%%%%%%%%%%%%%%
\section{The Integral Fractional Laplacian}\label{S:IntegralMethod}
%%%%%%%%%%%%%%%%%%%%%%%%%%%%%%%%%%%%%%%%%%%%%%%%%%%%%%%%%%%%%%%%%%%%%%%%%%%%%%%%%%%%

Here we consider the discretization of the integral definition of the fractional Laplacian \eqref{eq:Dir_frac}. %In what follows we assume that $\Omega$ is a bounded domain with Lipschitz boundary. 
In view of the definition \eqref{e:tildeHs} of the space $\Hs$,
and the fractional Poincar\'e inequality
\[
\| w \|_{L^2(\Omega)} \lesssim |w|_{H^s(\rn)} \quad \forall w \in \Hs,
\]
we may furnish $\Hs$ with the $H^s(\rn)$-seminorm. We also define the bilinear form $\llbracket\cdot,\cdot\rrbracket: \Hs \times \Hs \to \R$, 
\begin{equation}\label{sing-kernel}
\llbracket u, w \rrbracket := \frac{C(d,s)}{2} \iint_Q \frac{(u(x)-u(\othervar))(w(x)-w(\othervar))}{|x-\othervar|^{d+2s}} \diff \othervar \diff x,
\end{equation}
where $Q= (\Omega \times\rn) \cup (\rn \times \Omega)$ and $C(d,s)$ was defined in \eqref{eq:frac_lap}. We denote by $\interleave \cdot \interleave$ the norm that $\llbracket\cdot,\cdot\rrbracket$ induces, which is just a multiple of the $H^s(\rn)$-seminorm. The weak formulation of \eqref{eq:Dir_frac} is obtained upon multiplying \eqref{eq:frac_lap} by a test function $w\in\Hs$, integrating over $x\in \Omega$ and exploiting symmetry to make the difference $w(x)-w(\othervar)$ appear. With the functional setting we have just described at hand, this problem is formulated as follows: find $u \in \Hs$ such that
\begin{equation} \label{eq:cont}
  \llbracket u, w \rrbracket = \langle f, w\rangle  \quad \forall v \in \Hs.
\end{equation}
Applying the Lax-Milgram lemma immediately yields well-posedness of \eqref{eq:cont}.

Although the energy norm $\interleave \cdot \interleave$ involves integration on $\Omega\times\rn$, 
this norm can be localized. In fact, due to Hardy's inequality \cite{MR1993864,MR2085428}, the following equivalences hold \cite[Corollary 2.6]{AB}:
\[ \begin{aligned}
\| w \|_{H^s(\Omega)} & \le \interleave w \interleave  \lesssim \| w \|_{H^s(\Omega)}, & \text{ if } s \in (0,1/2) ,\\
| w |_{H^s(\Omega)} & \le  \interleave w \interleave  \lesssim | w |_{H^s(\Omega)}, & \text{ if } s \in (1/2,1).
\end{aligned} \]
When $s=1/2$, since Hardy's inequality fails, it is not possible to bound the $H^{\frac12}(\rn)$-seminorm in terms of the $H^{\frac12}(\Omega)$-norm for functions supported in $\overline \Omega$. However, for the purposes we pursue in this work, it suffices to notice that the estimate
\[
\interleave w \interleave \lesssim |w|_{H^{\frac12+\eps}(\Omega)} 
\]
holds for all $w \in \mathbb{H}^{\frac12+\eps}(\Omega)$.

From this discussion, it follows that the energy norm may be bounded in terms of fractional--order norms on $\Omega$. Thus, in order to estimate errors in the energy norm, we may bound errors within $\Omega$.

%-------------------------------------------------------------------------------------
\subsection{Regularity}\label{S:reg-integral}
%-------------------------------------------------------------------------------------
We now review some results regarding Sobolev regularity of solutions to problem \eqref{eq:cont} that are useful to deduce convergence rates of the finite element scheme proposed below. 
Regularity results for the fractional Laplacian have been recently obtained by Grubb \cite{Grubb} in terms of H\"ormander $\mu$-spaces \cite{Hormander}.  The work \cite{AB} has reinterpreted these in terms of standard Sobolev spaces. The following result, see \cite{BdPM,Grubb}, holds for domains with smooth boundaries, a condition that is too restrictive for a finite element analysis.

\begin{theorem}[smooth domains] \label{T:reg_grubb}
Let $s \in (0,1)$, $\Omega$ be a domain with $\pp\Omega \in C^\infty$, $f \in H^r(\Omega)$ for some $r\ge -s$, $u$ be the solution of \eqref{eq:Dir_frac} and $\gain = \min \{ s + r, 1/2 -\eps \}$, with $\eps > 0$ arbitrarily small. Then, $u \in \mathbb{H}^{s + \gain}(\Omega)$ and the following regularity estimate holds: 
\[
  \| u  \|_{\mathbb{H}^{s + \gain}(\Omega)} \lesssim \|f \|_{H^r(\Omega)},
\]
where the hidden constant depends on the domain $\Omega$, the dimension $d$, $s$ and $\gain$.
\end{theorem}

As a consequence of the previous result, we see that smoothness of the right hand side $f$ does not ensure that solutions are any smoother than $\mathbb{H}^{s + \frac12 - \eps}(\Omega)$; see also \cite{MR0214910}. We illustrate this phenomenon with the following example. 

\begin{example}[limited regularity]\label{ex:nonsmooth}
We follow \cite{Getoor,MR3447732} and consider $\Omega = B(0,1) \subset \rn$ and $f \equiv 1$. Then the solution to \eqref{eq:cont} is given by 
\begin{equation} \label{eq:getoor}
u(x) = \frac{\Gamma(\frac{d}{2})}{2^{2s} \Gamma(\frac{d+2s}{2})\Gamma(1+s)} ( 1- |x|^2)^s_+,
\end{equation}
where $t_+ =\max\{t,0\}$.
\end{example}

The lack of a lifting property for the solution to \eqref{eq:cont} can also be explained by the fact that the eigenfunctions of this operator have reduced regularity \cite{BdPM, Grubb_autovalores,RosOtonSerra2}. This is in stark contrast with
the spectral fractional Laplacian \eqref{def:second_frac}, discussed in Section \ref{S:reg-spectral},
whose eigenfunctions coincide with those of the Laplacian and thus are smooth functions if the boundary of the domain is regular enough.

On the other hand, H\"older regularity results for \eqref{eq:cont}
have been obtained in \cite{RosOtonSerra}. They give rise to Sobolev
estimates for solutions in terms of H\"older norms of the data, that
are valid for rough domains. More precisely, we have the following result; see \cite{AB}.

\begin{theorem}[Lipschitz domains] \label{T:regularity}
Let $s\in (0,1)$ and $\Omega$ be a Lipschitz domain satisfying the exterior ball condition. If $s\in (0,1/2)$,  let $f\in C^{\frac{1}{2}-s}(\overline\Omega)$; if $s = 1/2$, let $f \in L^\infty(\Omega)$; and if $s \in (1/2, 1)$, let $f\in C^{\beta}(\overline\Omega)$ for some $\beta>0$. Then, for every $\eps >0$, the solution $u$ of \eqref{eq:cont} belongs to $\mathbb{H}^{s+\frac{1}{2}-\eps}(\Omega)$, with
\[
  \|u\|_{\mathbb{H}^{s+\frac{1}{2}-\eps}(\Omega)} \lesssim  \frac1{\eps}  \|f\|_{\star},
\]
where $\| \cdot \|_\star$ denotes the $C^{\frac{1}{2}-s}(\overline\Omega)$, $L^\infty(\Omega)$ or $C^{\beta}(\overline\Omega)$, correspondingly to whether $s$ is smaller, equal or greater than $1/2$, and the hidden constant depends on the domain $\Omega$, the dimension $d$ and $s$.
\end{theorem}

In case $s>1/2$, the theorem above ensures that the solution $u$ belongs at least to $H^1_0(\Omega)$. 
It turns out that to prove the full
$H^{s+\frac{1}{2}-\eps}$-regularity, an intermediate step is to ensure
that the gradient of $u$ is actually an $L^2$-function. Following
\cite{BourgainBrezisMironescu}, this fact can be proved studying the
behavior of the fractional seminorms
$|\cdot|_{H^{1-\delta}({\Omega})}$, which usually blow up as $\delta \to 0$: 
\begin{equation} \label{eq:BBM}
\lim_{\delta \to 0} \ \delta |w|^2_{H^{1-\delta}({\Omega})} = C(d) |w|^2_{H^{1}({\Omega})} \quad \forall w \in L^2(\Omega).
\end{equation}
Therefore, the technique used in \cite{AB} to prove Theorem~\ref{T:regularity} consists of first proving that
the left-hand side of \eqref{eq:BBM}
remains bounded as $\delta \to 0$ for the solution $u$ of
\eqref{eq:cont}, whence $u \in H^1({\Omega})$, and
next analyzing the regularity of the gradient of $u$.

As already stated in \eqref{boundary-grubb}, the solution to \eqref{eq:cont} behaves like $\dist (x, \partial\Omega)^s$ for points $x$ close to the boundary $\partial\Omega$. This can clearly be seen in Example \ref{ex:nonsmooth} and explains the reduced regularity obtained in
Theorems \ref{T:reg_grubb} and \ref{T:regularity}. To capture this behavior, we develop estimates in fractional weighted  norms, where the weight is a power of the distance to the boundary. Following \cite{AB} we introduce the notation
\[
\delta(x,\othervar)=\min \big\{ \dist(x, \pp \Omega), \dist(\othervar, \pp \Omega) \big\},
\]
and, for $\ell = k + s$, with $k \in \mathbb{N}$ and $s \in (0,1)$, and $\kappa \ge 0$, we define the norm
\[
\| w \|_{H^{\ell}_\kappa (\Omega)}^2 := \| w \|_{H^k (\Omega)}^2 + \sum_{|\beta| = k } 
\iint_{\Omega\times\Omega} \frac{|D^\beta w(x)-D^\beta w(\othervar)|^2}{|x-\othervar|^{d+2s}} \, \delta(x,\othervar)^{2\kappa} \, \diff \othervar \, \diff x. 
\]
and the associated space
\begin{equation} \label{eq:weighted_sobolev}
H^\ell_\kappa (\Omega) :=  \left\{ w \in H^\ell(\Omega) \colon \| w \|_{H^\ell_\kappa (\Omega)} < \infty \right\} .
\end{equation}

Although we are interested in the case $\kappa \ge 0$, we recall that in the definition of weighted Sobolev spaces $H^k_\kappa(\Omega)$, with $k$ being a nonnegative integer, arbitrary powers of $\delta$ can be  considered \cite[Theorem 3.6]{Kufner}. On the other hand, global versions $H^\ell_\kappa(\rn)$ are defined integrating in the space $\rn$ and taking $\delta$ as before, but some restrictions must be taken into account to ensure their completeness. A sufficient condition is that the weight belongs to the Muckenhoupt class $A_2(\rn)$ \cite{Kilpelainen}. In this context, this implies that if $|\kappa|<1/2$ then the spaces $H^\ell_\kappa(\rn)$ are complete.

\begin{remark}[explicit solutions] \label{rem:explicit}
In radial domains, spaces like \eqref{eq:weighted_sobolev} have been used to characterize the mapping properties of the fractional Laplacian. In particular, when $\Omega$ is the unit ball, upon defining the weight $\omega(x) = 1 - |x|^2$, an explicit eigendecomposition of the operator $w \mapsto (-\Delta)^s (\omega^s w)$ is obtained in \cite{ABBM, DydaKuznetsovKwasnicki}.
The eigenfunctions are products of solid harmonic polynomials and radial Jacobi polynomials or, in one dimension, Gegenbauer polynomials. Mapping properties of the fractional Laplacian can thus be characterized in terms of weighted Sobolev spaces defined by means of expansions on these polynomials.
\end{remark}

The regularity in the weighted Sobolev spaces
\eqref{eq:weighted_sobolev} reads as follows \cite[Proposition 3.12]{AB}.

\begin{theorem}[weighted Sobolev estimate] \label{T:weighted_regularity}
Let $\Omega$ be a bounded, Lipschitz domain satisfying the exterior ball condition, $s\in(1/2,1)$, $f \in C^{1-s}(\overline\Omega)$ and $u$ be the solution of \eqref{eq:cont}. Then, for every $\eps>0$ we have $u \in H^{1+s-2\eps}_{1/2-\eps}(\Omega)$ and
\[
\|u\|_{H^{1+s-2\eps}_{1/2-\eps}(\Omega)} \lesssim \frac1\eps \| f \|_{C^{1-s}(\overline\Omega)},
\]
where the hidden constant depends on the domain $\Omega$, the dimension $d$ and $s$.
\end{theorem}

It is not straightforward to extend Theorem \ref{T:weighted_regularity} to $s\in (0,1/2]$ since, in this case, we cannot invoke Theorem \ref{T:regularity} to obtain that the solution $u$ belongs to $H^1(\Omega)$. Circumventing this would require to introduce a weight to obtain, for some $\kappa>0$, that $u  \in H^{1-\epsilon}_\kappa(\Omega)$ for some $\kappa>0$. However, for this is necessary to obtain a weighted version of \eqref{eq:BBM} which, to the best of our knowledge, is not available in the literature. 
In spite of this, the numerical experiments we have carried out using graded meshes and $s \in (0, 1/2]$ show the same order of convergence as for $s\in (1/2,1)$; see Table \ref{tab:ejemplo} below.

The proof of Theorem \ref{T:weighted_regularity} presented in \cite{AB} also shows that, if $f \in C^\beta(\overline\Omega)$ for some $\beta \in (0, 2-2s)$, $\ell < \min\{ \beta + 2s, \kappa + s + 1/2\}$, then we have
\begin{equation} \label{eq:weighted_estimate}
|u|_{H^\ell_\kappa(\Omega)} \lesssim \frac1{(\beta+\ell-2s)(1+2(\kappa-s-\ell))}
|f|_{C^\beta(\overline\Omega)},
\end{equation}
where the hidden constant depends only on $\Omega$ and the dimension $d$.
This shows that increasing the exponent $\kappa$ of the weight allows for the differentiability order $\ell$ to increase as well. In principle, there is no restriction on $\kappa$ above; however, 
in the next subsection we exploit this weighted regularity by introducing approximations on a family of graded meshes. There we show that the order of convergence (with respect to the number of degrees of freedom) is only incremented as long as $\kappa < 1/2$.

%-------------------------------------------------------------------------------------
\subsection{FEM: A Priori Error Analysis}\label{S:2-FEM-apriori}
%-------------------------------------------------------------------------------------

The numerical approximation of the solution to \eqref{eq:cont} presents an immediate difficulty: the kernel of the bilinear form \eqref{sing-kernel} is singular, and consequently, special care must be taken. For this reason, most existing approaches restrict themselves to the one dimensional case ($d=1$). Explorations in this direction can be found using finite elements \cite{DEliaGunzburger},
finite differences \cite{HuangOberman} and Nystr\"om methods \cite{ABBM}. For several dimensions the literature is rather scarce. A Monte Carlo algorithm that avoids dealing with the singular kernel was proposed in \cite{Kyprianou}.

Here we present a direct finite element approximation in arbitrary
dimensions. Following \S \ref{S:1-FEM-apriori}  we denote by $\T$ a
conforming and shape regular mesh
of $\Omega$, consisting of simplicial elements $K$ of diameter bounded
by $h$. In order to present a unified approach for the whole range $s \in (0,1) $ we only consider approximations of \eqref{eq:cont} by continuous functions. For $s <1/2$ it is possible to consider piecewise constants, but we do not explore this here.
With the finite element space $\U(\T)$ defined as in \eqref{eq:defofmathbbU},
the finite element approximation of \eqref{eq:cont} is then the unique solution to the problem: find $U \in \U(\T)$ such that
\begin{equation}
\label{eq:disc}
  \llbracket U, W \rrbracket = \langle f, W \rangle \quad \forall \, W \in \U(\T).
\end{equation}
From this formulation it immediately follows that $U$ is the projection (in the energy norm) of $u$ onto $\U(\T)$. Consequently, we have a C\'ea-like best approximation result
 \[
 \interleave u - U \interleave = \inf_{W \in \U(\T)} \interleave u - W \interleave.
 \]
Thus, in order to obtain a priori rates of convergence, it just
remains to bound the energy-norm distance between the discrete spaces
and the solution. One difficult aspect of dealing with fractional
seminorms is that they are not additive with respect to domain
decompositions. Nevertheless, it is possible to localize these norms
\cite{Faermann}: for all $w \in H^s(\Omega)$ we have
\[
|w|_{H^s(\Omega)}^2 \leq \sum_{K \in \T} \left[ \int_K \int_{S_K} \frac{|w (x) - w (\othervar)|^2}{|x-\othervar |^{d+2s}} \, \diff \othervar \diff x + \frac{C(d,\sigma)}{s h_K^{2s}} \| w \|^2_{L^2(K)} \right], 
\]
where $S_K$ is the patch associated with $K\in\T$ and $\sigma$ denotes the shape-regularity parameter of the mesh $\T$.
From this inequality it follows that, to obtain a priori error estimates, it suffices to compute interpolation errors over the set of
patches $\{K \times S_K \}_{K \in \T}$. The reduced regularity of solutions implies that we
need to resort to quasi-interpolation operators; we work with the
Scott-Zhang operator $\SZ$ \cite{ScottZhang}. Local stability and approximation properties of this operator were studied by Ciarlet Jr.  in \cite{Ciarlet}. 

\begin{proposition}[quasi-interpolation estimate] \label{prop:app_SZ}
Let $K \in \T$, $\max\{ 1/2, s \} <\ell \le 2$, $s \in (0,1)$, and $\SZ$ be the Scott-Zhang operator. 
If $w \in H^\ell (\Omega)$, then 
\[
 \int_K \int_{S_K} \frac{|(w-\SZ w) (x) - (w-\SZ w) (\othervar)|^2}{|x-\othervar |^{d+2s}} \,  \diff \othervar \, \diff x \lesssim
  h_K^{2\ell-2s} |w|_{H^\ell(S_K)}^2,
\]
where the hidden constant depends on $d$, $\sigma$, $\ell$ and blows
up as $s \uparrow 1$.
\end{proposition}

The interpolation estimate of Proposition \ref{prop:app_SZ} shows that, if the meshsize is sufficiently small, we deduce an a priori error bound in the energy norm.

\begin{theorem}[energy error estimate for quasi-uniform meshes] \label{T:conv_uniform}
Let $u$ denote the solution to \eqref{eq:cont} and denote by $U \in \U(\T)$ the solution of the discrete problem  \eqref{eq:disc}, computed over a mesh $\T$ consisting of elements with maximum diameter $h$. Under the hypotheses of Theorem \ref{T:regularity} we have
\[
  \interleave u - U \interleave \lesssim
   h^\frac12 |\log h| \|f\|_{\star},
\]
where the hidden constant depends on $\Omega$, $s$ and $\sigma$, and $\| \cdot \|_\star$ denotes the $C^{\frac{1}{2}-s}(\overline\Omega)$, $L^\infty(\Omega)$ or $C^{\beta}(\overline\Omega)$, correspondingly to whether $s$ is smaller, equal or greater than $1/2$.
\end{theorem}

This estimate hinges on the regularity of solutions provided by Theorem \ref{T:regularity}, and thus it depends on H\"older bounds for the data. We now turn our attention to obtaining a priori error estimates in the $L^2(\Omega)$-norm. Using Theorem \ref{T:reg_grubb}, an Aubin-Nitsche duality argument can be carried out. The proof of the following proposition follows the steps outlined in~\cite[Proposition 4.3]{BdPM}.

\begin{proposition}[$L^2$-error estimate]\label{prop:conv_L2}
Let $u$ denote the solution to \eqref{eq:cont} and denote by $U \in \U(\T)$ the solution of the discrete problem  \eqref{eq:disc}, computed over a mesh $\T$ consisting of elements with maximum diameter $h$. Under the hypotheses of Theorem \ref{T:reg_grubb} we have
\[
\| u - U \|_{L^2(\Omega)} \lesssim h^{\gain+\beta} \| f \|_{H^r(\Omega)} ,
\]
where $\gain = \min \{ s+r, 1/2-\eps \}$, $\beta = \min \{s, 1/2-\eps \},$ $\eps > 0 $ may be taken arbitrarily small and the hidden constant depends on $\Omega$, $s$, $d$, $\sigma$, $\gain$ and blows up when $\eps \to 0$.
\end{proposition}

Finally, for $s\in (1/2,1)$ and $d=2$, we take advantage of Theorem
\ref{T:weighted_regularity}, from which further information about the
boundary behavior of solutions is available. We
propose a standard procedure often utilized in connection with corner
singularities or boundary layers arising in convection-dominated problems. An increased rate of convergence is achieved by resorting to a priori adapted meshes.  To obtain interpolation estimates in the weighted fractional  Sobolev spaces defined by \eqref{eq:weighted_sobolev}, we introduce the following  Poincar\'e inequality \cite[Proposition 4.8]{AB}.

\begin{proposition}[weighted fractional Poincar\'e inequality] 
Let $s\in (0,1)$, $\kappa \in [0, s)$
and a domain $S$ which is star-shaped with respect to a ball. Then, for every $w \in H^s_\kappa (S)$, it holds
\[
\| w - \overline w \|_{L^2(S)} \lesssim d_S^{s-\kappa}|w|_{H^s_{\kappa}(S)},
\]
where $\overline w = \fint_S w \diff x$, $d_S = \diam (S)$ and the hidden constant depends on the chunkiness parameter of $S$ and the dimension $d$.
\end{proposition}

This inequality yields sharp quasi-interpolation estimates near the
boundary of the domain, where the weight $\delta$ involved in
\eqref{eq:weighted_sobolev} degenerates. Such bounds in turn
lead to
estimates for the Scott-Zhang operator in weighted fractional
spaces. We exploit them for two-dimensional problems ($d=2$) by constructing \emph{graded} meshes as in
\cite[Section 8.4]{Grisvard}. In addition to shape regularity, we assume that our meshes $\T$ satisfy the following property: there is a number $\mu\ge1$ such that given a meshsize parameter $h$ and $K\in\T$, we have
 \begin{equation} \label{eq:H}
 h_K \leq C(\sigma)
 \begin{dcases}
   h^\mu, & K \cap \partial \Omega \neq \emptyset, \\
   h \dist(K,\pp \Omega)^{(\mu-1)/\mu}, & K \cap \partial \Omega = \emptyset.
 \end{dcases}
\end{equation} 
where the constant $C(\sigma)$ depends only on the shape regularity constant $\sigma$ of the mesh $\T$. 
The parameter $\mu$ relates the meshsize $h$ to the number of degrees of freedom because (recall that $d=2$)
\[
 \#\T \approx
  \begin{dcases}
    h^{-2}, & \mu \leq 2, \\
    h^{-\mu}, & \mu > 2.
  \end{dcases}
\]
It is now necessary to relate the parameter $\mu$ with the exponent $\kappa$ of the weight $\delta$ in estimate \eqref{eq:weighted_estimate}. Increasing the parameter $\mu$ corresponds to raising $\kappa$ and thereby allowing an increase of the differentiability order $\ell$. However, if $\mu > 2$ this gain is compensated by a growth in the number of degrees of freedom. Following \cite{AB}, it turns out that the optimal parameter is $\mu=2$ and we have the following result.

\begin{theorem}[energy error estimates for graded meshes] \label{T:conv_graded}
Let $\Omega \subset \R^2$ and $U \in \U(\T)$ be the solution to \eqref{eq:disc}, computed over a mesh $\T$ that satisfies \eqref{eq:H} with $\mu = 2$. In the setting of Theorem \ref{T:weighted_regularity} we have
\[
  \interleave u - U \interleave  \lesssim  
  %\frac{C(s,\sigma)}{2s-1} 
  (\#\T)^{-\frac12} |\log (\#\T)|^\frac12 \, \|f\|_{C^{1-s}(\overline\Omega)}, 
\]
where the hidden constant depends on %$d$, 
$\sigma$ and blows up as $s \to 1/2$.
\end{theorem}

%---------------------------------------------------------------------------------
\subsection{Implementation} \label{sub:implementation}
%---------------------------------------------------------------------------------

Let us now discuss key details about the finite element
implementation of \eqref{eq:cont} for $d=2$. If $\{ \phi_\vero \}$ are the nodal piecewise linear basis functions of $\U(\T)$, defined as in \eqref{eq:defofmathbbU}, then the entries of the stiffness matrix $\A = (\A_{\vero\wero})$ are
\[
\A_{\vero\wero} = \llbracket \phi_\vero,\phi_\wero \rrbracket = \frac{C(d,s)}{2} \iint_Q \frac{(\phi_\vero(x) - \phi_\vero(\othervar))(\phi_\wero(x)-\phi_\wero(\othervar))}{|x-\othervar|^{2+2s}} \, \diff \othervar \, \diff x.
\]

Two numerical difficulties --- coping with integration
on unbounded domains and handling the \emph{non-integrable} singularity
of the kernel --- seem to discourage a direct finite element
approach. However, borrowing techniques from the boundary element
method \cite{SauterSchwab}, it is possible to compute accurately
the entries of the matrix $\A$. We next briefly outline the main steps of this
procedure. For full details, we refer to \cite{ABB}
where a finite element code to solve \eqref{eq:cont} is documented.

The integrals involved in the computation of $\A$ should be
carried over $\R^2$. For this reason it is convenient to consider a
ball $B$ containing $\Omega$ and such that the distance from
$\overline \Omega$ to $B^c$ is an arbitrary positive number.  This is
needed in order to avoid difficulties caused by lack of symmetry when
dealing with the integral over $\Omega^c$ when $\Omega$ is not a
ball. Together with $B$, we introduce an auxiliary triangulation ${
  \T_A}$ on $B\setminus\Omega$ such that the  complete triangulation
$\T_B$ over $B$ (that is $\T_B = \T \cup \T_A$)
remains admissible and shape-regular. 

We define, for $1\le \ell, m \le \# \T_B$ and $K_\ell,K_m \in \T_B$,
\begin{equation}
\begin{split}
I^{\vero,\wero}_{\ell,m} & := \int_{K_\ell}\int_{K_m} \frac{(\phi_\vero(x) - 
\phi_\vero(\othervar))(\phi_\wero(x)-\phi_\wero(\othervar))}{|x-\othervar|^{2+2s}} \, \diff \othervar \, \diff x, \\
J^{\vero,\wero}_\ell & := \int_{K_\ell}\int_{B^c} 
\frac{\phi_\vero(x)\phi_\wero(x)}{|x-\othervar|^{2+2s}} \, \diff\othervar \, \diff x,  \label{eq:integrales}
\end{split}
\end{equation}
whence we may write
\[ 
\A_{\vero\wero} =  \frac{C(d,s)}{2} \sum_{\ell = 1}^{\#\T_B} \left( 
\sum_{m=1}^{\#\T_B} I^{\vero,\wero}_{\ell,m} + 2 J^{\vero,\wero}_\ell \right).
\]
We reiterate that computing the integrals $I_{\ell, m}^{\vero,\wero}$ and
$J_\ell^{\vero,\wero}$ is challenging for different reasons: the former
involves a singular integrand if
$\overline{K}_\ell\cap\overline{K}_m\neq \emptyset$, while the latter
needs to be calculated in an unbounded domain.

We first tackle the computation of $I_{\ell,m}^{\vero,\wero}$ in
\eqref{eq:integrales}. If $\overline{K}_\ell$ and $\overline{K}_m$ do
not touch, then the integrand is a regular function and can be
integrated numerically in a standard fashion. On the other hand, if
$\overline K_\ell \cap \overline K_m \neq \emptyset$, then
$I_{\ell,m}^{\vero,\wero}$ bears some resemblances to  typical integrals
appearing in the boundary element method. Indeed, the quadrature rules
we employ are analogous to the ones presented in
\cite[Chapter 5]{SauterSchwab}. Basically, the scheme consists of the following
steps:
\begin{enumerate}[$\bullet$]
\item Consider parametrizations $\chi_\ell \colon \hat{K} \to K_\ell$
  and $\chi_m \colon \hat{K} \to K_m$ such that the edge or vertex shared
  by $K_\ell$ and $K_m$ is the image of the same edge/vertex in the reference element $\hat{K}$. If $K_\ell$ and $K_m$ coincide, simply use the same parametrization twice.
  
\item Decompose the integration domain into certain subsimplices and then utilize Duffy-type transformations to map these subdomains into the four-dimensional unit hypercube. 

\item Since the Jacobian of these Duffy transformations is regularizing, each of the integrals may be separated into two parts: a highly singular but explicitly integrable part and a smooth, numerically tractable part.
\end{enumerate}

The second difficulty lies in the calculation of $J^{\vero,\wero}_\ell$, namely, dealing with the unbounded domain $B^c$. We write
\[
J^{\vero,\wero}_\ell= \int_{K_\ell} \phi_\vero(x)\phi_\wero(x) \varrho(x) \, \diff x,
\qquad
\varrho(x):= \int_{B^c} \frac{1}{|x-\othervar|^{2+2s}} \, \diff \othervar,
\]
and realize that we need to accurately compute $\varrho(x)$ at
each quadrature point
$x \in K_\ell \cap \bar\Omega$. To do so, there are two properties we can take advantage of: the radiality of $\varrho$ and the fact that $\varrho$ is smooth up to the
boundary of $\Omega$ because, for $x \in \bar\Omega$ and $\othervar \in B^c$ $|x-\othervar|> \dist (\overline \Omega, B^c)>0.$ Therefore, the values of $\varrho$ at quadrature nodes can be precomputed with an arbitrary degree of
precision.

%-------------------------------------------------------------------------------------
\subsection{Numerical Experiments}\label{S:2-FEM-num-exp}
%-------------------------------------------------------------------------------------

We present the outcome of two experiments posed in $\Omega = B(0, 1) \subset \R^2$. 

%-------------------------------------------------------------------------------------
\subsubsection{Rate of Convergence in Energy Norm}\label{E1}
%-------------------------------------------------------------------------------------
Following Example \ref{ex:nonsmooth} we have that, if $f \equiv 1$, then the solution to \eqref{eq:cont} is given by \eqref{eq:getoor}. 
Table \ref{tab:ejemplo} shows computational rates of convergence in the energy norm for
several values of $s$, both for uniform and graded meshes. These
rates are in agreement with those predicted
by Theorems \ref{T:conv_uniform} and \ref{T:conv_graded}. Moreover,
we observe an increased order of convergence for $s\le 1/2$ and 
graded meshes which is not accounted for in Theorem \ref{T:conv_graded}.
\begin{table}[htbp]\small\centering
\begin{tabular}{|c| c| c| c| c| c| c| c| c| c|} \hline
Value of $s$ & 0.1 & 0.2 & 0.3 & 0.4 & 0.5 & 0.6 & 0.7 & 0.8 & 0.9 \\ \hline
Uniform meshes & 0.497 & 0.496 &  0.498  & 0.500 & 0.501  & 0.505 & 0.504  & 0.503 &  0.532 \\ \hline
Graded meshes & 1.066 & 1.040  & 1.019 & 1.002  & 1.066 & 1.051  & 0.990 & 0.985 & 0.977  \\  \hline   
\end{tabular}
\bigskip
\caption{Example \ref{E1}: Computational rates of convergence for
  \eqref{eq:cont} posed in the unit ball with right-hand side $f \equiv  1$. Errors are measured in the energy norm with respect to the
  meshsize parameter $h$. The second row corresponds to uniform meshes, while the third to graded meshes, with $\mu = 2$ in \eqref{eq:H}.}\label{tab:ejemplo}
\end{table}

%-------------------------------------------------------------------------------------
\subsubsection{Rate of Convergence in $L^2$-Norm}\label{E2}
%-------------------------------------------------------------------------------------
Remark \ref{rem:explicit} states that a family of explicit solutions
for \eqref{eq:cont} is available in the unit ball. A subclass of
solutions in that family may be expressed in terms of the Jacobi
polynomials $P_k^{(\alpha, \beta)}\colon [-1,1]\to\R$. We set
\[
f(x) = \left( \frac{\Gamma(3+s)}{2^{1-s}} \right)^2 P_2^{(s,0)}(2|x|^2 -1) ,
\]
so that the solution to \eqref{eq:cont} is given by \cite[Theorem 3]{DydaKuznetsovKwasnicki}
\[
u(x) = (1 - |x|^2)^s_+ \, P_2^{(s,0)}(2|x|^2 -1).
\]
We compute the orders of convergence in $L^2(\Omega)$ for  $s \in \{
0.25, 0.75 \}$;  according to Proposition \ref{prop:conv_L2}, it is
expected to have order of convergence $0.75$ for $s = 0.25$ and $1$
for $s=0.75$ with respect to the meshsize $h$. The results,
summarized in Table \ref{tab:resultados}, agree 
with the predicted rates of convergence.

\begin{table}[htbp]
\centering
\begin{tabular}{| c | c | c|}
\hline
$h$ & $s=0.25$ & $s=0.75$ \\ \hline
$0.0383$  & $0.0801$   & $0.01740$   \\ 
$0.0331$  & $0.0698$   & $0.01388$   \\ 
$0.0267$  & $0.0605$   & $0.01104$   \\ 
$0.0239$  & $0.0556$   & $0.00965$   \\ 
$0.0218$  & $0.0513$   & $0.00849$   \\ \hline
\end{tabular}
\bigskip
\caption{Example \ref{E2}:
Errors in the $L^2$-norm for $s = 0.25$ and $s = 0.75$. The estimated orders of convergence with respect to the meshsize are, respectively, $0.7669$ and $1.2337$.} \label{tab:resultados}
\end{table}

%-------------------------------------------------------------------------------------
\subsection{FEM: A Posteriori Error Analysis}\label{S:2-FEM-aposteriori}
%-------------------------------------------------------------------------------------

Since solutions of \eqref{eq:cont} have reduced regularity, and
assembling the stiffness matrix $\A$ entails a rather high
computational cost, it is of interest to devise suitable AFEMs.
We now present a posteriori error estimates of residual type
and ensuing AFEM; we follow \cite{NvPZ}.

We estimate the energy error $\interleave u-U \interleave$ in terms of the residual
$\calR := f - (-\Delta)^s U$ in $\Hsd$. To do so, we
need to address two important issues: localization of the norm in
$\Hsd$ and a practical computation of $\calR$.

To localize fractional norms we deviate from \cite{Faermann} and
perform a decomposition on stars $S_\vero = \supp \, \phi_\vero$, the
support of the basis functions $\phi_\vero$ associated with node $\vero$
and diameter $h_\vero$.
Exploiting the partition of
unity property $\sum_\vero \phi_\vero=1$, and Galerkin orthogonality
$\llbracket u- U,\phi_\vero \rrbracket=0$ for all $\vero \in\Omega$, we can write, for every $w \in \Hs$
\[
  \llbracket u - U, w \rrbracket = \langle \calR, w \rangle = \sum_\vero \langle \calR, w \phi_\vero \rangle
= \sum_\vero \langle \calR, (w- \bar w_\vero) \phi_\vero \rangle = \sum_\vero \langle (\calR-\bar \calR_\vero)\phi_\vero, w-\bar w_\vero\rangle,
\]
where $\bar w_\vero\in\R$ are weighted mean values computed as
$\bar w_\vero = 0$ provided $\vero\in\partial\Omega$ and, otherwise,
\[
\bar w_\vero := \frac{\langle w, \phi_\vero \rangle}{\langle \phi_\vero, 1  \rangle}
\quad\forall \, \vero \in\Omega.
\]
The values of $\bar \calR_\vero \in \R$ are yet to be chosen.
We see that $\calR = \sum_\vero (\calR-\bar \calR_\vero)\phi_\vero$ where each term
$(\calR-\bar \calR_\vero)\phi_\vero \in \mathbb{H}^{-s}(S_\vero)$
has support in $S_\vero$.
We have the
following two estimates for dual norms proved in \cite[Lemmas 1 and 2]{NvPZ}.

\begin{lemma}[localized upper bound of dual norm]\label{L:localized-upper}
Let $\calG\in \Hsd$ be decomposed as $\calG=\sum_\textup{\vero} g_\textup{\vero}$ with $g_\textup{\vero}\in \mathbb{H}^{-s}(S_\textup{\vero})$
vanishing outside $S_\textup{\vero}$. We then have for $0<s<1$
\[
\|\calG\|_{\Hsd}^2 \le (d+1) \sum_\textup{\vero} \|g_\textup{\vero}\|_{\mathbb{H}^{-s}(S_\textup{\vero})}^2.
\]
\end{lemma}

A key practical issue is then how to evaluate $\|g_\vero\|_{\mathbb{H}^{-s}(S_\vero)}$ for
$g_\vero=(\calR-\bar \calR_\vero)\phi_\vero$. For second order operators,
$\calR$ splits into an $L^2$-component in element
interiors (provided $f\in L^2(\Omega)$) and a singular component
supported on element boundaries. In contrast, the residual $\calR=f-(-\Delta)^s U$ 
does not have a singular component and its absolutely continuous part
is not always in $L^2(\Omega)$ for all $0<s<1$ no matter how smooth $f$ might be.
This is related to singularities of $(-\Delta)^s \phi_\vero(x)$ as $x$
tends to the skeleton of $S_\vero$ because $\phi_\vero$ is continuous, piecewise
linear. Using that \cite[Theorem XI.2.5]{MR618463}
\[
(-\Delta)^s: \widetilde W^t_p(\Omega) \rightarrow W^{t-2s}_p(\Omega), \qquad t \in \mathbb{R}, \quad p > 1,
\]
is a continuous pseudo-differential operator of
order $2s$, and $\phi_\vero \in \widetilde W^{1+\frac{1}{p}-\eps}_p(\Omega)$
for any $\eps>0$, we deduce
\[
(-\Delta)^s U \in L^p(\Omega) \qquad \frac{1}{p} > 2s-1.
\]
This motivates the following estimate, whose proof is given in
\cite[Lemma 2]{NvPZ}.

\begin{lemma}[upper bound of local dual norm]\label{L:upper-practical}
Let $g_\textup{\vero} \in L^p(S_\textup{\vero})$ satisfy $\int_{S_\textup{\vero}} g_\textup{\vero} \diff x=0$ for each
$\textup{\vero}\in\Omega$. If $0<s<1$ and $1\le p<\infty$ satisfies
$\frac{1}{p} < \frac{s}{d} + \frac{1}{2}$, then
\[
\|g_\textup{\vero}\|_{\mathbb{H}^{-s}(S_\textup{\vero})} \lesssim h_\textup{\vero}^{s+ \frac{d}{2}-\frac{d}{p}} \|g_\textup{\vero}\|_{L^p(S_\textup{\vero})}.
\]
\end{lemma}
However, to be able to apply Lemma \ref{L:upper-practical} the
Lebesgue exponent $p$ must satisfy
\[
2s - 1 < \frac{1}{p} < \frac{s}{d} + \frac{1}{2}.
\]
We note that for $s<\frac{3}{4}$ we can choose $p=2$ for any dimension
$d$. However, for $\frac{3}{4} \le s <1$ we need to take $1<p<2$. For
$d=1,2$ this condition is satisfied for any $s<1$, but for $d=3$ we
have the unfortunate constraint $s<\frac{9}{10}$.

We can now choose $\bar \calR_\vero$ as
\[
  \bar{\calR}_\vero :=  \frac{ \int_{S_\vero} \calR \phi_\vero \diff x }{\int_{S_\vero} \phi_\vero \diff x}, \quad \vero \in \Omega
\]
and $\bar \calR_\vero = 0$ otherwise,
so that the local contributions satisfy $\int_{S_\vero} (\calR - \bar\calR_\vero)\phi_\vero \diff x= 0$ and we can then apply the bound of Lemma \ref{L:upper-practical}.
The following upper a posteriori error estimate is derived in
\cite[Theorem 1]{NvPZ}.

\begin{theorem}[upper a posteriori bound]\label{T:upper-bound}
  Let $f\in L^p(\Omega)$ and $1<p<\infty, 0<s<1$ satisfy the restriction
  $2s-1 < \frac{1}{p} < \frac{s}{d}+\frac{1}{2}$, then
  \[
  \|u-U\|_{\Hsd}^2 \lesssim \sum_\textup{\vero}
  h_\textup{\vero}^{2\big(s+ \frac{d}{2}-\frac{d}{p}\big)}  \|(\calR-\bar \calR_\textup{\vero})\phi_\textup{\vero}\|_{L^p(S_\textup{\vero})}^2.
  \]
\end{theorem}

This error analysis has two pitfalls. The first one, alluded to
earlier, is a restriction on $s$ for $d>2$. The second one is the actual
computation of $(-\Delta)^s \phi_\vero(x)$ for $d>1$, which is
problematic due to its singular behavior as $x$ tends to the skeleton
on $\T$. This is doable for $d=1$ and we refer to \cite[Section 8]{NvPZ} for
details and numerical experiments. This topic is obviously open for improvement.

%-------------------------------------------------------------------------------------
 \subsection{Extensions and Applications}
%-------------------------------------------------------------------------------------

We conclude the discussion by mentioning extensions and applications of this approach:
\begin{enumerate}[$\bullet$]
\item \emph{Eigenvalue problems}: The eigenvalue problem for the integral fractional Laplacian arises, for example, in the study of fractional quantum mechanics \cite{Laskin}. As already mentioned, a major difference between the spectral fractional Laplacian and the integral one is that, 
%while for the first operator eigenspaces coincide with those of the Dirichlet Laplacian --- and thus consist of smooth functions under standard regularity assumptions on the domain --- 
for the second one, eigenfunctions have reduced regularity.  In \cite{BdPM}, conforming finite element approximations were analyzed and it was shown that the Babu\v{s}ka-Osborn theory \cite{BO91} holds in this context. Regularity results for the eigenfunctions are derived under the assumption that the domain is Lipschitz and satisfies the exterior ball condition. Numerical evidence on the non-convex domain $\Omega = (-1,1)^2\setminus[0,1)^2$ indicates that the first eigenfunction is as regular as the first one on any smooth domain. This is in contrast with the Laplacian.
\item \emph{Time dependent problems}: In \cite{ABB2}, problem \eqref{eq:CaPUTOivbp} with $\gamma \in (0,2]$ and the integral definition of $\Laps$ was considered. Regularity of solutions was studied and a discrete scheme was proposed and analyzed. The method is based on a standard Galerkin finite element approximation in space, as described here, while in time a convolution quadrature approach was used \cite{JinLazarovZhou,Lubich}.
\item \emph{Non-homogeneous Dirichlet conditions}: An interpretation of a non-homogeneous Dirichlet condition $g$ for the integral fractional Laplacian is given by using \eqref{eq:frac_lap} upon extension by $g$ over $\Omega^c$. In \cite{ABH} a mixed method for this problem is proposed; it is based on the weak enforcement of the Dirichlet condition and the incorporation of a certain non-local normal derivative as a Lagrange multiplier. This non-local derivative is interpreted as a non-local flux between $\Omega$ and $\Omega^c$ \cite{DROV17}. 
\item \emph{Non-local models for interface problems}: Consider two
  materials with permittivities/diffusivities of opposite sign, and
  separated by an interface with a corner. Strong singularities may
  appear in the classical (local) models derived from electromagnetics
  theory. In fact, the problem under consideration is of Fredholm type if and only if the quotient between  the value of permittivities/diffusivities taken from both sides of the interface lies outside a so-called critical interval, which always contains the value $-1$. In \cite{BC} a non-local interaction model for the materials is proposed. Numerical evidence indicates that the non-local model may reduce the critical interval and that solutions are more stable than for the local problem.
\end{enumerate} 

% !TEX root = three-methods.tex
%%%%%%%%%%%%%%%%%%%%%%%%%%%%%%%%%%%%%%%%%%%%%%%%%%%%%%%%%%%%%%%%%%%%%%%%%%%%%%%%%%%%
\section{Dunford-Taylor Approach for Spectral and Integral Laplacians}\label{S:DumfordTaylor}
%%%%%%%%%%%%%%%%%%%%%%%%%%%%%%%%%%%%%%%%%%%%%%%%%%%%%%%%%%%%%%%%%%%%%%%%%%%%%%%%%%%%

%As already mentioned in the introduction, we present 
In this section we present an alternative approach to the ones
developed in the previous sections.
It relies on the Dunford-Taylor representation \eqref{e:balak_inv}
$$
u = (-\Delta)^{-s}f = \frac{\sin(s \pi)}{\pi } \int_{0}^\infty \mu^{-s} (\mu-\Delta)^{-1}f \diff\mu
$$
for the spectral fractional Laplacian \eqref{def:second_frac}. 
For the integral fractional Laplacian \eqref{eq:Fourier}, instead,
it hinges on the equivalent representation \eqref{e:rep_int} of \eqref{e:weak_int_lap}:
\begin{equation}\label{e:equiv_bil2}
\begin{split}
&\int_{\mathbb R^d} |\xi|^s \mathscr{F}^{-1}(\tilde u)  |\xi|^s \overline{\mathscr{F}^{-1}(\tilde w)} \diff \xi \\
&\qquad \qquad = 
\frac{2\sin(s\pi)}{\pi} \int_0^\infty \mu^{1-2s} \int_{\Omega} \left((-\Delta)(I-\mu^2 \Delta)^{-1}\tilde u(x) \right) w(x) \diff x \diff \mu.
\end{split}
\end{equation}
In \eqref{e:equiv_bil2}, the operators $-\Delta$ and $I-\mu^2\Delta$
are defined over $\mathbb R^d$, something to be made precise in
Theorem \ref{T:dunf:int:equv}.

In each case, the proposed method is proved to be efficient on general Lipschitz domains $\Omega \subset \mathbb R^d$.
They rely on sinc quadratures and on finite element approximations of the resulting integrands at each quadrature points.
While \eqref{e:balak_inv} allows for a direct approximation of the solution, the approximation of \eqref{e:rep_int} leads to a non-conforming method where the action of the stiffness matrix on a vector is approximated.

We recall that the functional spaces $\Hs$ are defined in
\eqref{e:tildeHs} for $s\in [0,3/2)$. We now extend the definition
$\Hs = H^s(\Omega) \cap H^1_0(\Omega)$ for $s \in (1,2]$.
%they are equivalent to $(L^2(\Omega),H^1_0(\Omega))_s$ when $s \in [0,1]$ and to $H^s(\Omega) \cap H^1_0(\Omega)$ when $s \in (1,3/2)$.

%---------------------------------------------------------------------------------
\subsection{Spectral Laplacian}\label{ss:dunford-spectral}
%---------------------------------------------------------------------------------

We follow \cite{BP1,BP2} and describe a method based on the Balakrishnan representation \eqref{e:balak_inv}.
In order to simplify the notation, we set $L:=-\Delta: \mathcal{D}(L)
\rightarrow L^2(\Omega)$ and define the domain of $L^r$, for $r \in \R$,
to be
$$
\mathcal{D}(L^r):= \{ v \in L^2(\Omega) \ : \ L^rv \in L^2(\Omega) \};
$$
this is a Banach space equipped with the norm 
\begin{equation*}
%\label{e:normDA}
\| v\|_{\mathcal{D}(L^r)} := \| L^r v \|_{L^2(\Omega)}.
\end{equation*}

We also define the solution operator $T: H^{-1}(\Omega) \rightarrow H^1_0(\Omega)$ by $Tf:=v$, where for $F \in H^{-1}(\Omega)$, $v\in H^1_0(\Omega)$ is the unique solution of 
$$
\int_\Omega \nabla v \cdot \nabla w \diff x = F(v), \qquad \forall w \in H^1_0(\Omega).
$$
This definition directly implies that $\mathcal{D}(L) = \textrm{range}(T|_{L^2(\Omega)})$.

%----------------------------------------------------------------------------------
\subsubsection{Finite Element Discretization}\label{ss:dunford-spectral-fem}
%----------------------------------------------------------------------------------

For simplicity, we assume that the domain $\Omega$ is polytopal so that it can be partitioned into a conforming subdivision $\T$.
We recall that $\mathbb U(\T) \subset H^1_0(\Omega)$  stands for the subspace of globally continuous piecewise linear polynomials with respect to $\T$; see Section~\ref{S:1-FEM-apriori}. 
We denote by $\Pi_{\T}$ the $L^2(\Omega)$-orthogonal projection onto $\mathbb U(\T)$ and by  $L_\T:\mathbb U(\T) \rightarrow \mathbb U(\T)$ the finite element approximation of $L$, i.e.,  for $V \in \mathbb U(\T)$, $L_\T V \in \mathbb U(\T)$ solves
$$
\int_\Omega (L_\T V) W\diff x = \int_\Omega \nabla V \cdot \nabla W\diff x, \qquad \forall\, W \in \mathbb U(\T).
$$
We finally denote by $T_\T$ the inverse of $L_\T$, the finite element
solution operator, and by $h$ the maximum diameter of elements in $\T$. 

With these notations, we are in the position to introduce the finite element approximation $U \in \mathbb U(\T)$ of $u$ in \eqref{e:dunf-inv}:
\begin{equation}\label{e:dunf:discrete_spectral}
U := \frac{\sin(s\pi)}{\pi} \int_0^\infty \mu^{-s} (\mu +L_\T)^{-1} \Pi_{\T}f \diff \mu.
\end{equation}

The efficiency of the approximation of $u$ by $U$ depends on the efficiency of the finite element solver $(T_\T \Pi_\T) f$ (i.e. for the standard Laplacian), which is dictated by the regularity of $Tf$.
This regularity aspect has been intensively discussed in the literature \cite{bramblebacuta,MR3085047,dauge2,Joch99,kellogg,others}.
In this exposition, we make the following general assumption.

\begin{definition}[elliptic regularity]
We say that $T$ satisfies a pick-up regularity of index $0<\alpha\leq 1$ on $\Omega$ if for $0\leq r \leq \alpha$,  the operator 
$T$ is an isomorphism from $\mathbb H^{-1+r}(\Omega)$ to $\mathbb H^{1+r}(\Omega)$.
\end{definition}

Notice that $\alpha =1$ when $\Omega$ is convex,
whence this definition extends \eqref{Omega_regular} to general Lipschitz domains.

Assuming a pick-up regularity of index $\alpha$,  for any $r \in [0,1]$,  we have
$$
\| Tw - T_\T \Pi_\T w \|_{\mathbb H^r(\Omega)} \lesssim
h^{2\alpha_*} \|Tw\|_{\mathbb H^{\alpha+1}(\Omega)} \lesssim
h^{2\alpha_*} \|w\|_{\mathbb H^{\alpha-1}(\Omega)}
$$
where $\alpha_* := \frac 1 2 \big(\alpha + \min(1-r,\alpha)\big)$.
The proof of the above estimate is classical and is based on a duality argument (Nitsche's trick); see e.g. \cite[Lemma~6.1]{BP2}. 
Notice that $\alpha_* < \alpha$ when $1-\alpha<r$, i.e, the error is
measured with regularity index too large to take full advantage of the pick-up regularity in the duality argument.

We expect that approximation \eqref{e:dunf:discrete_spectral} of the
fractional Laplacian problem delivers the \emph{same} rate of convergence
\begin{equation}\label{e:Dunf:heuristic_space}
\|u-U\|_{\mathbb{H}^r(\Omega)} \lesssim
h^{2\alpha_*} \|u\|_{\mathbb{H}^{1+\alpha}(\Omega)}
\lesssim h^{2\alpha_*} \|f\|_{\mathbb{H}^{1+\alpha-2s}(\Omega)},
\end{equation}
but the function $U$ in \eqref{e:dunf:discrete_spectral} is well defined provided $f \in L^2(\Omega)$, i.e. $1+\alpha-2s\ge0$. 

Before describing the finite element approximation result, we make the following comments. 
The solution $u = L^{-s}f$ belongs to $\mathcal{D}(L^{(1+\alpha)/2})$ provided that $f \in \mathcal D(L^{(1+\alpha)/2-s})$.
Hence, estimates such as \eqref{e:Dunf:heuristic_space} rely on the characterization of $\mathcal {D}(L^{r/2})$ for $r \in (0,1+\alpha]$.
For $0\leq r \leq 1$, the spaces $ \mathcal{D}(L^{r/2})$ and $\mathbb H^r(\Omega)$ are equivalent, as they are both scale spaces which coincide at $r = 0 $ and $r =1$. 
Furthermore,  assuming an elliptic regularity pick-up of index $0<\alpha\leq 1$, this characterization extends up to $1+\alpha$ \cite[Theorem~6.4 and Remark~4.2]{BP2}.

 \begin{theorem}[finite element approximation]\label{t:dunf:spec_fem}
Assume that $T$ satisfies a pick-up regularity of index $\alpha \in (0,1]$ on $\Omega$.
Given $r\in [0,1]$ with $r \leq 2s$, set $\gamma :=\max\{r+2\alpha_*-2s,0\}$
and $\alpha_* := \frac12 \big( \alpha + \min (1-r,\alpha) \big)$.
If $f \in \mathbb H^{\delta}(\Omega)$ for $\delta \geq \gamma$,
then
$$
\|u-U \|_{\mathbb H^{r}(\Omega)} \leq C_{h} h^{2\alpha_*} \| f \|_{\mathbb H^{\delta}(\Omega)},
$$
where
$$
C_{h} \lesssim \left\lbrace
\begin{array}{ll}
\log(2/h), \qquad & \textrm{when } \delta = \gamma \quad \textrm{and} \quad r+2\alpha_* \geq 2s,\\% \quad r+\alpha_* \not = 1/2, \\
1, \qquad & \textrm{when } \delta > \gamma \quad \textrm{and} \quad r+2\alpha_* \geq 2s, \\
1, \qquad & \textrm{when } \delta = 0 \quad \textrm{and} \quad 2s > r+2\alpha_*.
\end{array}
\right.
$$
\end{theorem}

%--------------------------------------------------------------------------------
\subsection{Sinc Quadrature}\label{ss:dunf:spect:sinc}
%--------------------------------------------------------------------------------

It remains to put in place a sinc quadrature \cite{lundbowers} to approximate the integral in \eqref{e:dunf:discrete_spectral}. 
We use the change of variable $\mu = e^y$ so that
$$
U = \frac{\sin(s\pi)}{\pi} \int_{-\infty}^\infty e^{(1-s)y} (e^y +L_\T)^{-1} \Pi_{\T}f \diff y.
$$
Given $k>0$, we set 
$$
N_+ := \left\lceil   \frac{\pi^2}{4s k^2} \right\rceil, \qquad  
N_- := \left\lceil   \frac{\pi^2}{4(1-s) k^2} \right\rceil, \qquad \text{and}\qquad y_\ell := k\ell,
$$
and define the sinc quadrature approximation of $U$  by
\begin{equation}\label{e:dunf:spec_full}
U^k := \frac{\sin(s\pi)}{\pi} k \sum_{\ell=-N_-}^{N_+} e^{(1-s)y_\ell} (e^{y_\ell} +L_\T)^{-1} \Pi_{\T}f.
\end{equation}
The sinc quadrature consists of uniformly distributed quadrature
points in the $y$ variable, and the choice of $N_+$ and $N_-$ makes it more robust
with respect to $s$.

The decay when $|z| \to +\infty$ and holomorphic properties of the
integrand $z^{-s} (z - L)^{-1}$ in the Dunford-Taylor representation
\eqref{e:dunf-inv} guarantee the \emph{exponential convergence} of the
sinc quadrature \cite[Theorem~7.1]{BP2}.

\begin{theorem}[sinc quadrature]\label{t:dunf:int:sinc}
For  $r \in [0,1]$, we have
$$
\| U - U^k \|_{\mathbb H^r(\Omega)} \lesssim e^{-\pi^2/(2k)} \| f \|_{\mathbb H^r(\Omega)}. 
$$
\end{theorem}

To compare with Theorem \ref{TH:fl_error_estimates}, we take $r=s$ and
assume that $\Omega$ is convex, which allows for any $\alpha$ in  $(0,1]$.
We choose a number of sinc quadrature points $N_+ \approx N_- \approx
\log(1/h)$ so that sinc quadrature and finite element errors are balanced. 
Therefore, Theorems~\ref{t:dunf:spec_fem} and \ref{t:dunf:int:sinc} yield for the Dunford-Taylor method
$$
\| u - U^k \|_{\mathbb H^s(\Omega)} \lesssim (\#\T)^{-2\alpha_*/d} \| f \|_{\mathbb H^\sigma(\Omega)},
$$
where $\sigma := \max(2\alpha_*-s,s)$ and $(\#\T)^{-1/d} \approx h$ for
quasi-uniform subdivisions, provided we discard logarithmic terms. 
In contrast, the error estimate of Theorem
\ref{TH:fl_error_estimates} for the extension method reads, again
discarding log terms,
\begin{equation*}
\| u - U \|_{\Hs} \lesssim (\# \T_\Y)^{-1/(d+1)} \|f \|_{\mathbb{H}^{1-s}(\Omega)}
\end{equation*}
and was derived with pick-up regularity $\alpha=1$.
We first observe the
presence of the exponent $d+1$, which makes the preceding error
estimate suboptimal. This can be cured with geometric grading in the
extended variable and $hp$-methodology. Section \ref{sub:ext_NOS} and
\cite{BMNOSS:17,MFSV:17} show that this approach yields the following error estimate 
\begin{equation*}
\| u - U \|_{\Hs} \lesssim (\# \T)^{-1/d} \|f \|_{\mathbb{H}^{1-s}(\Omega)},
\end{equation*}
with $\# \T$ degrees of freedom, after discarding logarithmic terms.
This estimate exhibits quasi--optimal linear order for the regularity
$f \in\mathbb{H}^{1-s}(\Omega)$.
We also see that the Dunford-Taylor method
possesses the optimal rate of convergence $2\alpha_*=2-s > 1$ allowed
by polynomial interpolation theory for smoother datum $f\in \mathbb{H}^\sigma(\Omega)$:
$\sigma = 2(1-s)$ when $s \leq 2/3$ and $\sigma =s$ when $s >2/3$.
We may also wonder what regularity of $f$ would lead to the same
linear order of convergence as the extension method, that requires
$f\in\mathbb{H}^{1-s}(\Omega)$. We argue as follows: if $s\le\frac12$,
then $f \in \Ws$; otherwise, if $s>\frac12$,
then $f\in \mathbb{H}^s(\Omega)$. We thus realize that the regularity of
$f$ is the same for $s\le\frac12$ but it is stronger for $s>\frac12$.
  
It is worth mentioning that the Dunford-Taylor algorithm seems
advantageous in a multi-processor context as it appears to exhibit
good strong and weak scaling properties.
The former consists of increasing the number of processors for a fixed
number of degrees of freedom, while for the latter, the number of
degrees of freedom per processor is kept constant when increasing the
problem size.
We refer to  \cite{CPE:CPE4216} for a comparison of different methods.

\begin{remark}[implementation]
The method based on \eqref{e:dunf:spec_full} requires $N_++N_-+1$ \emph{independent} standard Laplacian finite element solves for each quadrature points $y_\ell$:
$$
V^\ell \in \mathbb U(\T): \quad e^{y_\ell} \int_\Omega V^\ell  ~W\diff x + \int_\Omega \nabla V^\ell \cdot \nabla W \diff x = \int_\Omega f W \diff x \quad \forall W\in \mathbb U(\T),
$$
which are then aggregated to yield $U^k$:
$$
U^k = \frac{\sin(s\pi)}{\pi} k \sum_{\ell=-N_-}^{N_+} e^{(1-s)y_\ell} V^\ell.
$$
Implementation of this algorithm starting from a finite element solver
for the Poisson problem is straightforward.
Numerical illustrations matching the predicted convergence rates
of Theorems \ref{t:dunf:spec_fem} and \ref{t:dunf:int:sinc} are provided in \cite{BP1}.
\end{remark}

%----------------------------------------------------------------------------------
\subsubsection{Extensions}
%----------------------------------------------------------------------------------

We now discuss several extensions.

\begin{enumerate}[$\bullet$]
\item \emph{Symmetric operators and other boundary conditions.}
The operator $L=-\Delta$ can be replaced by any symmetric elliptic operators as long as 
the associated bilinear form $(u,w) \mapsto \int_\Omega Lu~w$  remains coercive and bounded in $H^1_0(\Omega)$.  
Different boundary conditions can be considered similarly as well.
However, it is worth pointing out that the characterization of $\mathcal{D}(L^r)$ depends on the boundary condition and must be established.

As an illustration, Figure~\ref{f:dunf_LB} depicts the approximations using parametric surface finite element \cite{Dziuk} of 
the solution to
\begin{equation}\label{e:LB}
(-\Delta_\Gamma)^{s} u =  1 \quad \textrm{on} \quad  \Gamma, \qquad u=0 \quad \textrm{on } \partial \Gamma,
\end{equation}
where $\Delta_\Gamma$ is the surface Laplacian on $\Gamma \subset \mathbb R^3$, either the side boundary of a cylinder or given by 
\begin{equation}\label{e:dunf:gamma}
\begin{split}
\Gamma := \bigl\{ (x_1+2\sin(x_3),x_2+2\cos(x_3)  ,10x_3) \in \mathbb
R^3  :  (x_1,x_2,x_3) \in \mathbb S_2,  x_3 \geq 0  \bigr\}.
\end{split}
\end{equation}
\begin{figure}[ht!]
\includegraphics[width=0.35\textwidth]{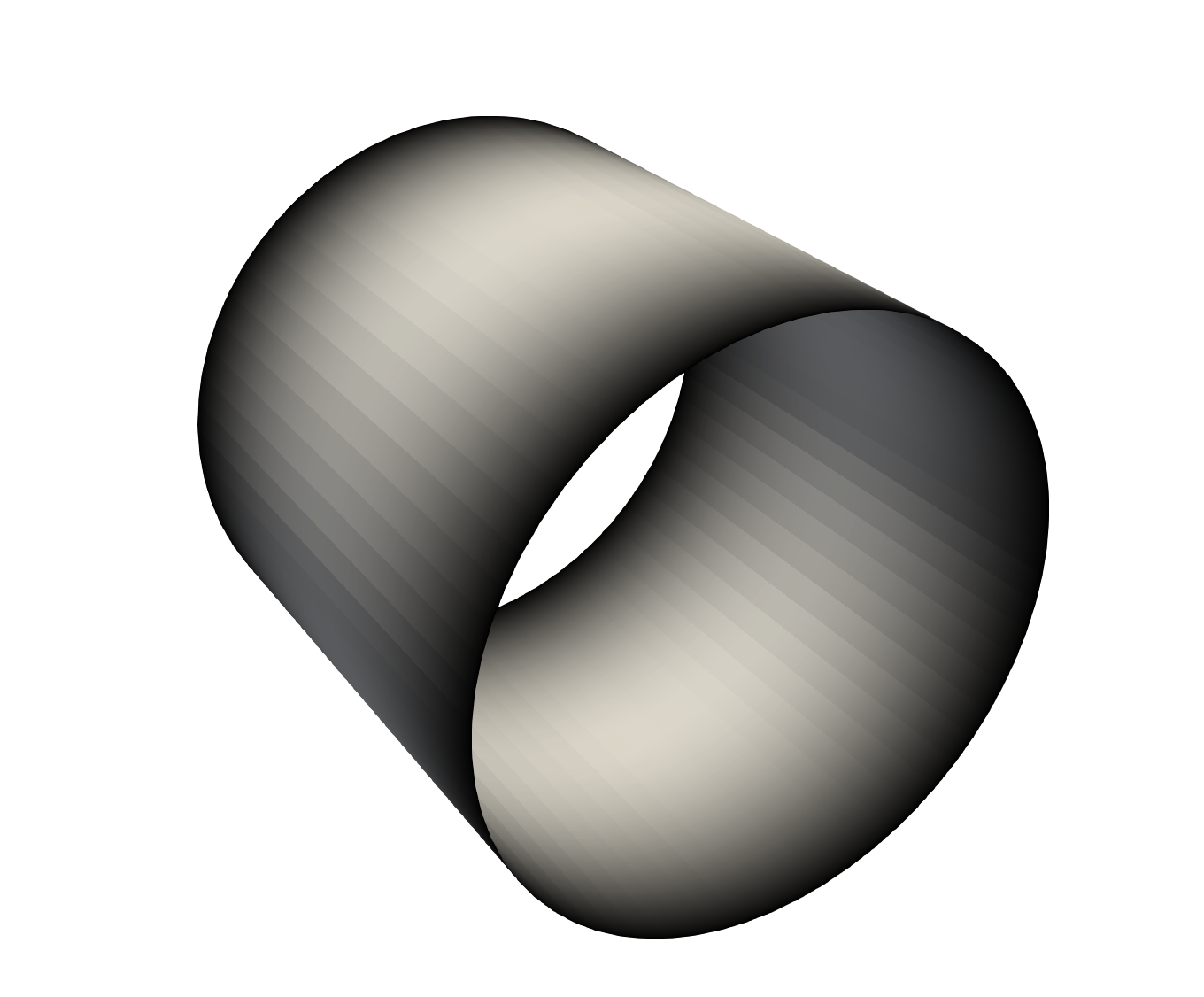}
\includegraphics[width=0.45\textwidth]{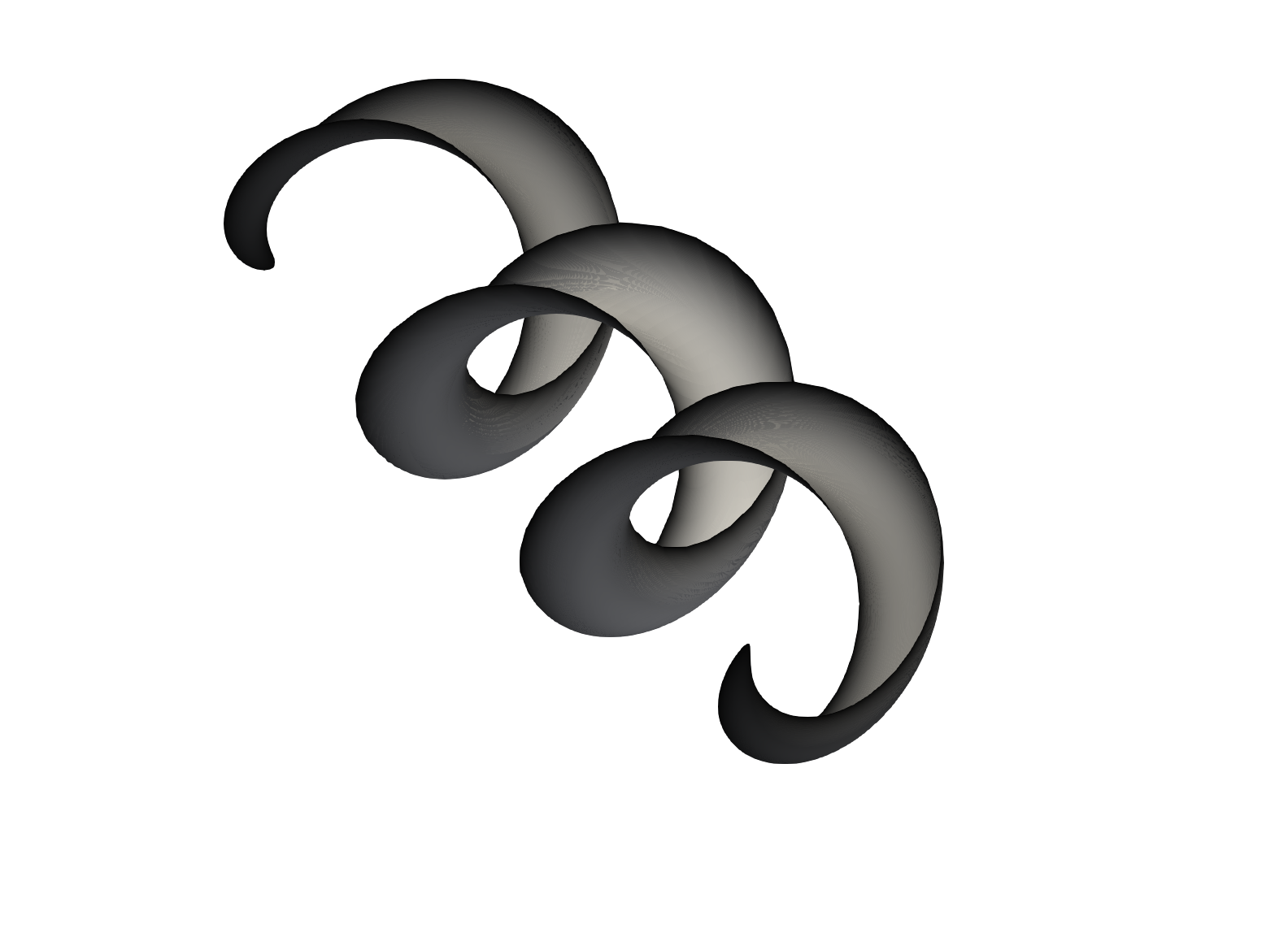}
\caption{Dunford-Taylor Method: Numerical approximation of the solution to the spectral Laplacian problem \eqref{e:LB} on an hypersurface (darker = smaller values;  lighter = larger values).
(Left) $s=0.8$ and $\Gamma$ is the side boundary of a cylinder of radius 1 and height 2; (Right) $s=0.5$ and $\Gamma$ is given by \eqref{e:dunf:gamma}.
}\label{f:dunf_LB}
\end{figure}

\item \emph{Regularly accretive operators.}  The class of operators $L$ can be extended further to a subclass of non-symmetric operators. 
They are the unbounded operators associated with coercive and bounded sesquilinear forms in $H^1_0(\Omega)$ (regularly accretive operators \cite{kato1961}). 
In this case, fractional powers cannot be defined using a spectral
decomposition as in \eqref{def:second_frac} but rather directly using
the Dunford-Taylor representation \eqref{e:dunf-inv} and the
Balakrishnan formula \eqref{e:balak_inv}, which remain valid. 
The bottleneck is the characterization of the functional spaces $\mathcal{D}(L^{r/2})$ in terms of Sobolev regularity.
It turns out that for $-1<r<1$, we have that $\mathcal{D}(L^{r/2})$ is the same as for the symmetric operator \cite{kato1960}
$$
\mathcal{D}(L^{r/2}) = \mathcal{D}((L+L^*)^{r/2}) = \mathbb H^{r}(\Omega),
$$
where $L^*$ denotes the adjoint of $L$.
This characterization does not generally hold for $r=1$ (Kato square root problem \cite{kato1961}).
McKintosh \cite{mcintosh90} proved that $\mathcal{D}(L^{1/2}) = \mathbb H^1(\Omega) = H^1_0(\Omega)$ for  sesquilinear forms of the type
$$
(\psi,\varphi) \mapsto  \int_\Omega \left(B \nabla \psi \cdot \nabla \overline \varphi + \beta_1 \cdot \nabla \psi \overline \varphi +  \psi \beta_2 \cdot \nabla \overline \varphi + c \psi \overline \varphi \right)\diff x,
$$
where $B \in L^\infty(\Omega, \GL(\rn))$, $\beta_1, \beta_2 \in L^\infty(\Omega,\mathbb R^d)$ and $c \in L^\infty(\Omega)$ are such that the form is coercive and bounded.
This characterization is extended in \cite[Theorem~6.4]{BP2} up to $r=1+\alpha$ so that similar convergence estimates to those in Theorems~\ref{t:dunf:spec_fem} and~\ref{t:dunf:int:sinc} are established. 
To illustrate the method for non-symmetric operators, we consider the following example
$$
\left(-\Delta  + \begin{pmatrix}1\\1 \end{pmatrix} \cdot \nabla\right)^{s} u = 1, \qquad \textrm{in }\Omega, \qquad u=0 \quad \textrm{on }  \partial \Omega,
$$
where $\Omega=(-1,1)^2 \setminus [0,1)^2$ is a L-shaped domain.
Figure~\ref{f:dunf-nonsym} reports the fully discrete approximations given by \eqref{e:dunf:spec_full} for $s=0.2,0.5,0.8$.
\begin{figure}[ht!]
\includegraphics[width=0.45\textwidth]{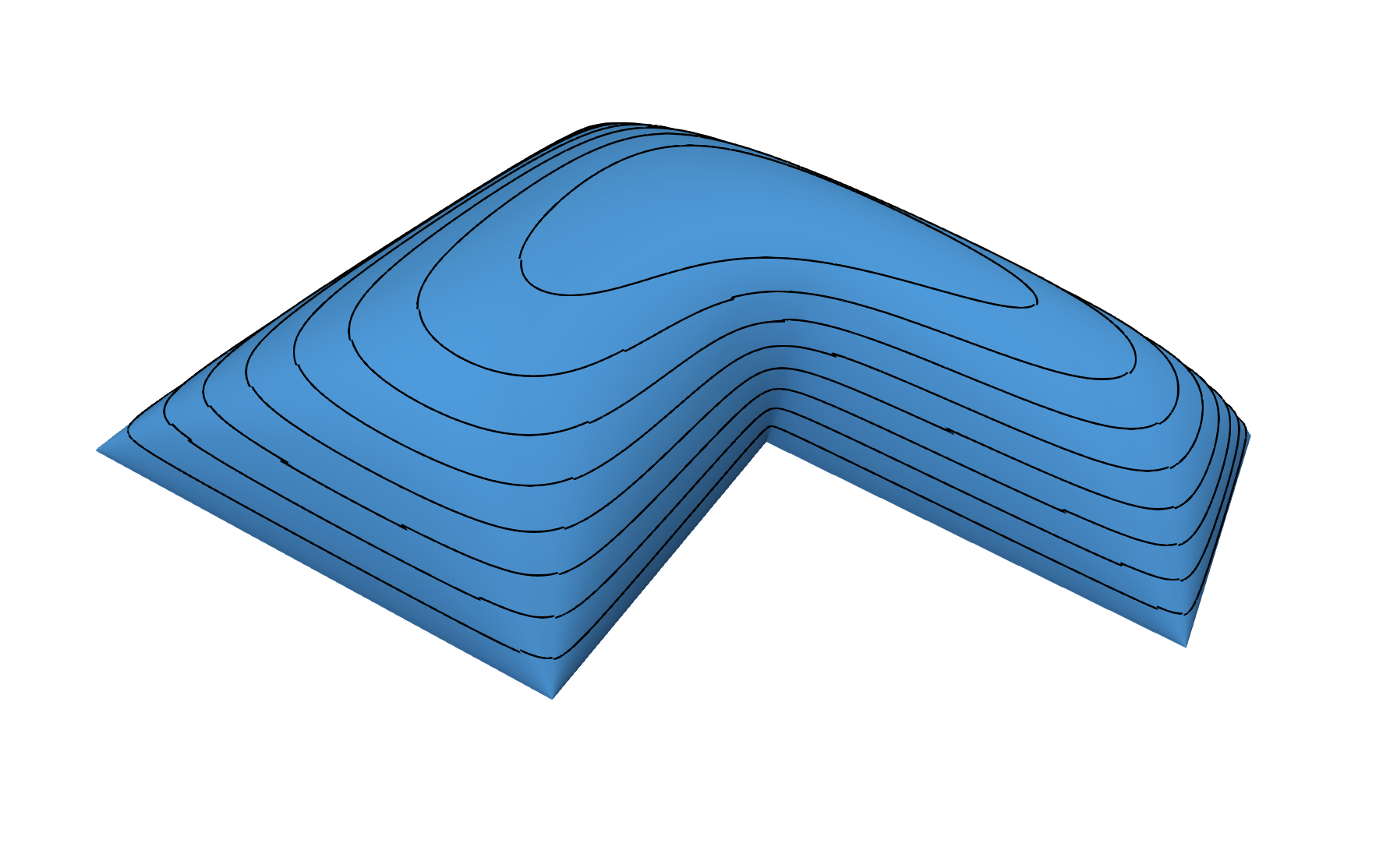}
\includegraphics[width=0.45\textwidth]{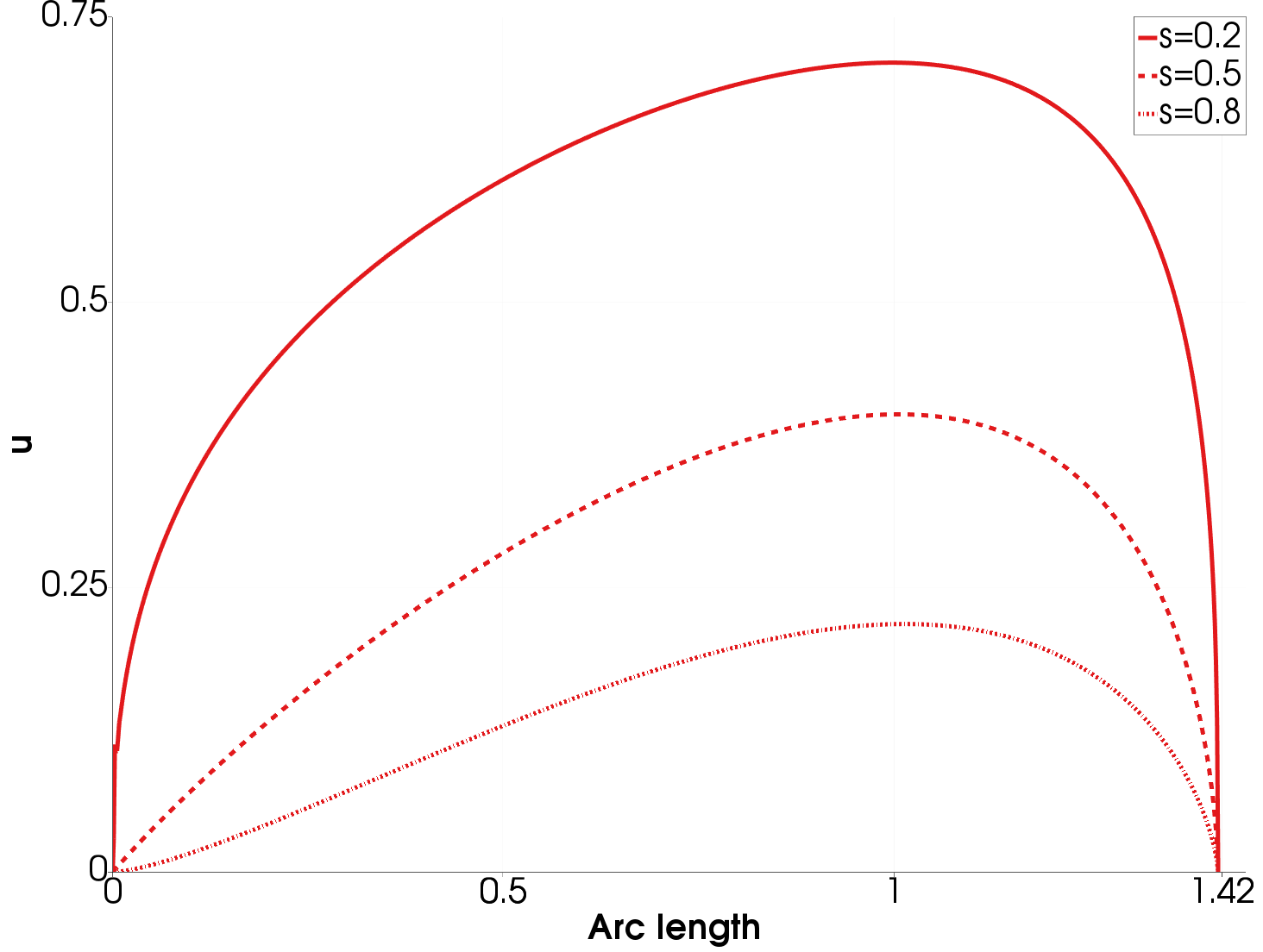}
\caption{Dunford-Taylor Method: Approximation of the solution to fractional advection - diffusion problem on a L-shaped domain. 
(Left) Solution with isolines for $s=0.5$. (Right) Plots of $u$ for $s=0.2,0.5,0.8$ on the segment from the corner opposite to the re-entrant corner (of coordinate (-1,1)) to the re-entrant corner (of coordinate (0,0)). It appears that the boundary layer intensity (but not its width) at the re-entrant corner depends on the power fraction $s$.}\label{f:dunf-nonsym}
\end{figure}

\item \emph{Space and time fractional diffusion.}
In \cite{BLP1,BLP2} the space-time fractional problem 
$$
\partial_t^\gamma u + (-\Delta)^s u = f, \qquad u(0) =u_0, 
$$
is studied, where $\partial_t^\gamma$ denotes the so-called Caputo derivative of order $\gamma \in (0,1]$; see \eqref{eq:caputo}.
The solution of the space-time fractional problem is given by \cite{SY11}
$$
u(t) = e_{\gamma,1}(-t^\gamma L^s) u_0 + \int_0^t \xi^{\gamma-1} e_{\gamma,\gamma}(-\xi^\gamma L^s) f(t-\xi)\diff\xi,
$$
where, again in this case, a Dunford-Taylor representation can be used to write
$$
e_{\gamma,\mu}(-t^\gamma L^s) = \frac{1}{2\pi i} \int_\contour e_{\gamma,\mu}(-t^\gamma z^s)(z-L)^{-1}\diff z
$$
and $e_{\gamma,\mu}$ defined on $\mathbb C$ is the Mittag-Leffler function.
Because of the presence of $e_{\gamma,\mu}$, the contour $\contour$ cannot be deformed  onto the negative real axis anymore, which prevents a representation like in \eqref{e:balak_inv}.
Instead, the sinc quadrature is performed directly  on a hyperbolic parametrization of $\contour$ in the complex plane.
Nevertheless, we obtain the error estimates
\begin{equation*}
 \| u(t) - U^k(t) \|_{L^2(\Omega)} 
  \lesssim \left(t^{-\gamma \alpha_*/s} h^{2\alpha_*} + t^{-\gamma} e^{-c/k} \right) \| u_0 \|_{L^2(\Omega)},
\end{equation*}
when $f=0$ and where $\alpha_*:= \min(\alpha,s^-)$ with $s^-$ denoting any number strictly smaller than $s$.
Here $c$ is a constant independent of $h$ and $k$.
We refer the reader to \cite{BLP1,BLP2} for estimates measuring the error in higher norms or for improved results (in the singularity when $t\to 0$) when $u_0$ is smoother. 
Note that the representation used does not need a time-stepping method for the initial value problem.

When $u_0=0$ instead,  a graded (a-priori known and depending on $\gamma$) mesh in time towards $t=0$ is put forward.
A midpoint quadrature scheme (second order) in time for a total of $\mathcal N \log(\mathcal N)$ time steps yields the error estimates
\begin{equation*}
\begin{split}
 \| u(t) - U^{k,\mathcal N} (t) \|_{L^2(\Omega)} \lesssim  & ~ t^{(1-\alpha_*/s)\gamma} h^{2\alpha_*} \|f \|_{L^2( (0,t)\times \Omega)}  \\
 &+ \max\{t^\gamma, t^{3/2+\gamma} \} \mathcal N^{-2} \| f\|_{H^2(0,t;L^2(\Omega))} \\
 & + \log(\mathcal N) e^{-c/k} \|f \|_{L^\infty(0,t;L^2(\Omega))}.
\end{split}
\end{equation*}
Notice that the method exhibits \emph{second order} convergence rate (up to a logarithmic term) with respect to the number of time intervals. 
Again, we refer to \cite{BLP2,BLP1} for more details as well as additional estimates when measuring the error in higher order norms.
\end{enumerate}

%----------------------------------------------------------------------------------
\subsection{Integral Laplacian}\label{ss:dunford-integral}
%----------------------------------------------------------------------------------

The strategy used for the spectral Laplacian in the previous section cannot be used for the integral Laplacian.
In fact, formulas like \eqref{e:dunf-inv} are not well defined (the integral Laplacian is not strictly positive).

Instead, we rely on the following equivalent representation of the bilinear form \eqref{e:equiv_bil2} in the weak formulation \eqref{e:weak_int_lap}.
Recall that fractional order Sobolev spaces in $\R^d$
are defined and normed by
$$
H^r(\R^d) = \left\lbrace w \in L^2(\R^d) : \ \| w \|_{H^r(\R^d)} = \left(\int_{\R^d} (1+|\xi|^2)^{r/2} | \mathscr{F}(w)(\xi)|^2 \diff \xi \right)^{1/2} < \infty \right\rbrace
$$
for $r>0$,
and that the notation $\tilde w$ stands for the zero extension of $w$ outside
$\Omega$ so that $w\in\mathbb{H}^r(\Omega)$ if and only if $\tilde
w\in H^r(\R^d)$ for $r\in [0,3/2)$; see definition \eqref{e:tildeHs}.

\begin{theorem}[equivalent representation]\label{T:dunf:int:equv}
Let $s \in (0,1)$ and $0\leq r \leq s$. For $\eta \in H^{r+s}(\R^d)$ and $\theta \in H^{s-r}(\R^d)$,
\begin{equation*}
\begin{split}
 &\int_{\R^d} |\xi|^{r+s} \mathscr{F}(\eta)(\xi) \, |\xi|^{s-r}  \overline{\mathscr{F}(\theta)(\xi)} \diff \xi  \\
 & \qquad \qquad \qquad = \frac{2\sin(s\pi)}{\pi} \int_0^\infty \mu^{1-2s} \int_{\R^d} \left((-\Delta)(I-\mu^2 \Delta)^{-1}\eta \right) \theta \diff x \diff\mu.
 \end{split}
 \end{equation*}
\end{theorem}
To prove the above theorem, it suffices to note that using Parseval's theorem
$$
 \int_{\R^d} \left((-\Delta)(I-\mu^2 \Delta)^{-1}\eta \right) \theta \diff x = \int_{\R^d} \frac{|\xi|^2}{1+\mu^2 |\xi|^2} \mathscr{F}(\eta)(\xi) \overline{\mathscr{F}(\theta)(\xi)}\diff \xi.
$$
and use the change of variable $t = \mu | \xi|$ together with the relation 
$$
\int_0^\infty \frac{t^{1-2s}}{1+t^2}\diff t = \frac{\pi}{2\sin(\pi s)}.
$$
For more details, we refer to \cite[Theorem~4.1]{BLP3}.

In order to make the above representation more amenable to numerical methods, for $\psi \in L^2(\R^d)$, we define  $v(\psi,\mu) := v(\mu) \in H^1(\R^d)$ to be the solution to
\begin{equation}\label{e:dunf:int_v}
\int_{\R^d} v(\mu)\phi \diff x+ \mu^2 \int_{\R^d} \nabla v(\mu) \cdot \nabla \phi \diff x= - \int_{\R^d} \psi \phi \diff x, \qquad \forall \phi \in H^1(\R^d).
\end{equation}
Using this notation along with definition \eqref{e:tildeHs} of $\Hs$,
we realize that for $\eta, \theta \in \Hs$ with $s\in (0,1)$, we have
$$
\int_{\R^d} |\xi|^{s} \mathscr{F}(\tilde \eta)(\xi) \, |\xi|^{s}
\overline{\mathscr{F}(\tilde \theta)(\xi)} \diff \xi =
\frac{2\sin(s\pi)}{\pi} \int_0^\infty
\mu^{-1-2s} \Big( \int_\Omega \big( \eta + v(\tilde \eta,\mu)\big) \theta \diff x\Big) \diff \mu;
$$
note that $v(\tilde \eta,\mu)$ does not vanish outside $\Omega$.
This prompts the definition
 \begin{equation}\label{e:a_int_lap}
 a(\eta,\theta) :=  \frac{2\sin(s\pi)}{\pi} \int_0^\infty \mu^{-1-2s} \Big( \int_\Omega \big( \eta + v(\tilde \eta,\mu)\big) \theta \diff x\Big) \diff \mu,
 \end{equation}
 for $\eta, \theta \in \Hs$. 
 The above representation is the starting point of the proposed numerical method. The solution $u \in \Hs$ of the fractional Laplacian satisfies
\begin{equation}\label{e:int_lap_weak_2}
a(u,w) = \langle f , w \rangle \qquad \forall w \in \Hs.
\end{equation}

We discuss in Section~\ref{ss:dunf:strang} a Strang's type argument to
assess the discretization error  from the consistency errors generated
by the approximation of $a(\cdot,\cdot)$ using sinc quadratures,
domain truncations, and finite element discretizations.

%------------------------------------------------------------------------------------
\subsubsection{Sinc Quadrature}
%------------------------------------------------------------------------------------

We proceed as in Section~\ref{ss:dunf:spect:sinc} for the spectral
Laplacian and use the change of variable $\mu = e^{-\frac 1 2 y}$
to arrive at
$$
a(\eta,\theta) = \frac{\sin(s\pi)}{\pi} \int_{-\infty}^\infty e^{sy} \left(\int_\Omega \big( \eta + v(\tilde \eta,\mu(y))\big) \theta \diff x\right) \diff y.
$$
Given a quadrature spacing $k>0$, $N_+$ and $N_-$ two positive integers, the sinc quadrature approximation of $a(\cdot,\cdot)$ is given by
\begin{equation}\label{e:dunf_int_sinc}
a^k(\eta,\theta) := \frac{\sin(s\pi)}{\pi} k \sum_{\ell = -N_-}^{N_+} e^{sy_\ell}
\int_\Omega \big( \eta + v(\tilde \eta,\mu(y_\ell))\big) \theta \diff x .
\end{equation}
Notice that we only emphasize the dependency in $k$ in the
approximation of $a(\cdot,\cdot)$ as we will select $N_+$ and $N_-$ as a function of $k$. 

The consistency error between $a^k(\cdot,\cdot)$ and
$a(\cdot,\cdot)$ is described in the following result.
We simply note that, as for the spectral Laplacian discussed in
Section~\ref{ss:dunford-spectral}, the proof of Theorem
\ref{t:dunf:int:quad} is given in \cite[Theorem~5.1 and Remark~5.1]{BLP3}
and relies on the holomorphic property and decay as $\mu \to \infty$
of the integrand in \eqref{e:a_int_lap}.

\begin{theorem}[quadrature consistency]\label{t:dunf:int:quad}
Given $\theta \in \Hs$ and $\eta \in \mathbb H^\delta(\Omega)$ with
$s<\delta  \leq \min(2-s,\sigma)$, where $\sigma$ stands for any number strictly smaller than $3/2$. 
Set $N_+ := \left\lceil \frac{\pi^2}{2k^2(\delta-s)}\right\rceil$ and $N_- := \left\lceil \frac{\pi^2}{4sk^2}\right\rceil$. Then, we have
$$
|a(\eta,\theta) - a^k(\eta,\theta)| \lesssim \max\left(\frac 1 {\delta-s}, \frac 1 s\right) e^{-\pi^2/(4k)} \| \eta \|_{\mathbb H^\delta(\Omega)} \| \theta \|_{\Hs}.
$$
\end{theorem}

%---------------------------------------------------------------------------------
\subsubsection{Truncated Problems}\label{ss:dunf:trunc}
%---------------------------------------------------------------------------------

The sinc approximation of $a(\cdot,\cdot)$ defined by \eqref{e:dunf_int_sinc} requires the computation of $v(\tilde\eta,\mu(y_\ell))$ for each quadrature point $y_\ell$ (here $\eta \in \mathbb H^\beta(\Omega)$ for some $s<\beta<3/2$ is fixed).
This necessitates, according to \eqref{e:dunf:int_v}, the approximations of $(I-\mu(y_\ell)^2\Delta)^{-1}$ on $\R^d$.
The proposed method relies on truncations of  this infinite domain problem and uses standard finite elements on the resulting bounded domains. 
As we shall see, the truncated domain diameter must  depend on the quadrature point $y_\ell$.

We let $B$ be a convex bounded domain containing $\Omega$ and the origin of $\R^d$.
Without loss of generality, we assume that the diameter of $B$ is 1. 
For a truncation parameter $M$, we define the dilated domains
\begin{equation*}%\label{e:BM} %JPB
B^M(\mu) := \left\lbrace
\begin{array}{ll}
\{ y = (1+\mu(1+M)x \ : x \in B \}, &\qquad \mu \geq 1, \\
\{ y = (2+M)x \ : x \in B \}, &\qquad \mu < 1,
\end{array}
\right.
\end{equation*}
and for $\psi \in L^2(\R^d)$, the associated functions $v^M(\mu) := v^M(\psi,\mu) \in H^1_0(B^M(\mu))$ satisfying
\begin{equation}\label{e:dunf:int_v_M}
\begin{split}
\int_{B^M(\mu)} v^M(\mu)w \diff x &+ \mu^2 \int_{B^M(\mu)} \nabla v^M(\mu) \cdot \nabla w \diff x \\
&\qquad = - \int_{B^M(\mu)} \psi w \diff x \qquad \forall w \in H^1_0(B^M(\mu));
\end{split}
\end{equation}
compare with \eqref{e:dunf:int_v}.
The exponential decay of the function $v(\tilde \eta,\mu)$ yields 
$$
\| v(\tilde \eta,\mu) - v^M(\tilde \eta,\mu) \|_{L^2(B^M(\mu))} \lesssim e^{-\max(1,\mu)cM/\mu} \| \eta \|_{L^2(\Omega)},
$$
where $c$ is a constant independent of $M$ and $\mu$ (see \cite[Lemma 6.1]{BLP3}).
As a consequence, the truncation consistency in using 
$$
a^{k,M}(\eta,\theta) := \frac{\sin(s\pi)}{\pi} k \sum_{\ell = -N_-}^{N_+} e^{sy_\ell} \int_\Omega (\eta + v^M(\tilde \eta,\mu(y_\ell))) \theta \diff x
$$
instead of $a^k(\eta,\theta)$ decays exponentially fast as a function
of $M$ \cite[Theorem~6.2]{BLP3}.

\begin{theorem}[truncation consistency]
For $M$ sufficiently large, there is positive constant $c$ independent of $M$ and $k$ such that
for all $\eta,\theta \in L^2(\Omega)$
$$
|a^k(\eta,\theta) - a^{k,M}(\eta,\theta)| \lesssim e^{-cM} \| \eta \|_{L^2(\Omega)} \| \theta \|_{L^2(\Omega)}.
$$
\end{theorem}

%---------------------------------------------------------------------------------
\subsubsection{Finite Element Discretization}
%---------------------------------------------------------------------------------

We now turn our attention to the finite element approximation of $v^M(\tilde \eta,\mu)$ defined by \eqref{e:dunf:int_v_M}.
For simplicity, we assume that the domain $\Omega$ is polytopal
so that it can be partitioned into a conforming subdivision $\T$ with
elements of maximum diameter $h$
as in Section \ref{ss:dunford-spectral-fem}.
Generic constants below may depend on the shape regularity and
quasi-uniformity constants of $\T$ without mention of it.

We need two subspaces of globally continuous piecewise linear
polynomials. The first one, $\mathbb U(\T) \subset H^1_0(\Omega)$,
is defined in \eqref{eq:defofmathbbU} relative to the partition
$\T$. The second subspace, denoted $\mathbb U(\T^M(\mu))$, has a
similar definition but relative to the subdivision $\T^M(\mu)$ of
$B^M(\mu(y_\ell))$. We impose that the partitions $\T^M(\mu)$ match in
$\Omega$, which implies that restrictions of functions in
$\mathbb U(\T^M(\mu))$ are continuous piecewise linears over $\T$. We
refer to \cite{BLP3} for details on the constructions of such
partitions, which is the bottleneck of the proposed method.

We now define for any $\psi \in L^2(\R^d)$ the finite element approximation
$V^M(\mu) = V^M(\psi,\mu) \in \mathbb U(\T^M(\mu))$ of the function
$v^M(\mu)$ to be
\begin{equation*}
\begin{split}
\int_{B^M(\mu)} V^M(\mu)W \diff x & + \mu^2 \int_{B^M(\mu)} \nabla V^M(\mu) \cdot \nabla W  \diff x \\
& \qquad \qquad = - \int_{B^M(\mu)} \psi W  \diff x \qquad \forall \, W \in \mathbb U(\T^M(\mu)).
\end{split}
\end{equation*}
The fully discrete approximation of the bilinear form $a(\cdot,\cdot)$ then reads
$$
a_\T^{k,M}(\Xi,\Theta) := \frac{\sin(s\pi)}{\pi} k \sum_{\ell =-
  N_-}^{N_+} e^{sy_\ell} \int_\Omega \big( \Xi + V^M(\widetilde
\Xi,\mu(y_\ell))\big) \Theta \diff x
\quad\forall \, \Xi, \Theta \in \mathbb U(\T).
$$
Note that $V^M(\widetilde \Xi,\mu)$ is piecewise linear
over $\T$ and the sum $\Xi + V^M(\widetilde \Xi,\mu)$ is easy to perform.
The consistency error between $a^{k,M}(\cdot,\cdot)$ and
$a^{k,M}_\T(\cdot,\cdot)$ is given next; see \cite[Theorem~7.6]{BLP3} for a proof.

\begin{theorem}[finite element consistency]
If $\beta \in (s,3/2)$ and $\alpha \in (0,\min(s,1/2))$, then the
following estimate is valid
$$
|a^{k,M}(\Xi,\Theta) - a^{k,M}_\T(\Xi,\Theta) | \lesssim |\log h| h^{\beta+\alpha-s} \|  \Xi \|_{\mathbb H^\beta(\Omega)} \| \Theta \|_{\mathbb H^{s+\alpha}(\Omega)},
$$
for all $\Xi,\Theta \in \mathbb U(\T)$.
\end{theorem}

%-----------------------------------------------------------------------------------
\subsubsection{Strang's Lemma}\label{ss:dunf:strang}
%-----------------------------------------------------------------------------------

In addition to the three consistency estimates described above, Strang's Lemma requires the $\mathbb U(\T)$-ellipticity of the fully discrete form $a_\T^{k,M}(\cdot,\cdot)$.
To show this, \cite{BLP3} invokes Theorem~\ref{t:dunf:int:quad}
(quadrature consistency) with $\delta :=\min\{2-s,\beta\}$ and an inverse estimate to write for $\Xi, \Theta \in \mathbb U(\T)$
$$
| a^k(\Xi,\Theta) - a(\Xi,\Theta) | \lesssim e^{-\pi^2/(4k)} h^{s-\delta} \| \Xi \|_{\mathbb H^s(\Omega)} \| \Theta \|_{\mathbb H^s(\Omega)}.
$$
This, together with the monotonicity property
$$
a^{k,M}_\T(\Xi,\Theta) \geq a^k(\Xi,\Theta),
$$
and the coercivity of the exact bilinear form $a(\cdot,\cdot)$ in
$\mathbb H^s(\Omega)$, yields the $\mathbb U(\T)$-ellipticity of
$a^{k,M}_\T(\cdot,\cdot)$ provided
\begin{equation}\label{e:dunf:cond}
e^{-\pi^2/(4k)} h^{s-\delta} \leq c
\end{equation}
for an explicit constant $c$. 

The fully discrete approximation $U^{k,M} \in \mathbb U(\T)$ of $u$ satisfying  \eqref{e:int_lap_weak_2} is given by
$$
a_\T^{k,M}(U^{k,M},W) = \int_\Omega f W \diff x\qquad \forall \, W  \in \mathbb U(\T).
$$
To measure the discrepancy between $u$ and $U^{k,M}$ in $\mathbb{H}^s(\Omega)$,
we assume the additional regularity $u \in \mathbb H^\beta(\Omega)$
for some $\beta \in (s,3/2)$. The expected
regularity of $u$, solution to the integral fractional Laplacian, is discussed
in Theorems~\ref{T:reg_grubb} and~\ref{T:regularity}.
The theorem below, proved in \cite[Theorem~7.8]{BLP3}), guarantees
that the proposed method delivers an \emph{optimal} rate of
convergence (up to a logarithmic factor).

\begin{theorem}[error estimate]\label{t:dunf:error_int}
Assume that \eqref{e:dunf:cond} holds and that $u \in \mathbb H^\beta(\Omega)$ for some $\beta \in (s,3/2)$.
Then there is a constant $c$ independent of $h$, $M$, and $k$ such that
$$
\| u - U^{k,M} \|_{\Hs} \lesssim \left(e^{-\pi^2/(4k)} 
 + e^{-cM} + |\log h| h^{\beta-s} \right) \| u \|_{\mathbb H^\beta(\Omega)}.
$$
\end{theorem}
We also refer to \cite[Theorem~7.8]{BLP3} for further discussions on mesh generation, matrix representation of the fully discrete scheme, and a preconditioned iterative method.

Notice that the error estimate stated in
Theorem~\ref{t:dunf:error_int} for quasi-uniform meshes is of order
about $h^{1/2}$ and is similar to
the one derived for the integral method in Theorem~\ref{T:conv_uniform}.
To see this, we choose $\beta = s+\frac 1 2-\varepsilon$ with
$\epsilon>0$ arbitrary,
which is consistent with the regularity of $u$ guaranteed by
Theorem~\ref{T:regularity}, along with $M \approx \log(1/h)$ and
$N_+ \approx N_- \approx |\log h|^2$ to balance the three sources of errors. 
It is worth mentioning that using graded meshes for $d=2$,
Theorem~\ref{T:conv_graded} states an optimal linear rate of convergence (up to a
logarithmic factor) provided $s\in (1/2,1)$.
Whether such a strategy applies to the Dunford-Taylor method
remains open.

%---------------------------------------------------------------------------------
\subsubsection{Numerical Experiment}
%---------------------------------------------------------------------------------

\begin{figure}[ht!]
\includegraphics[width=0.55\textwidth]{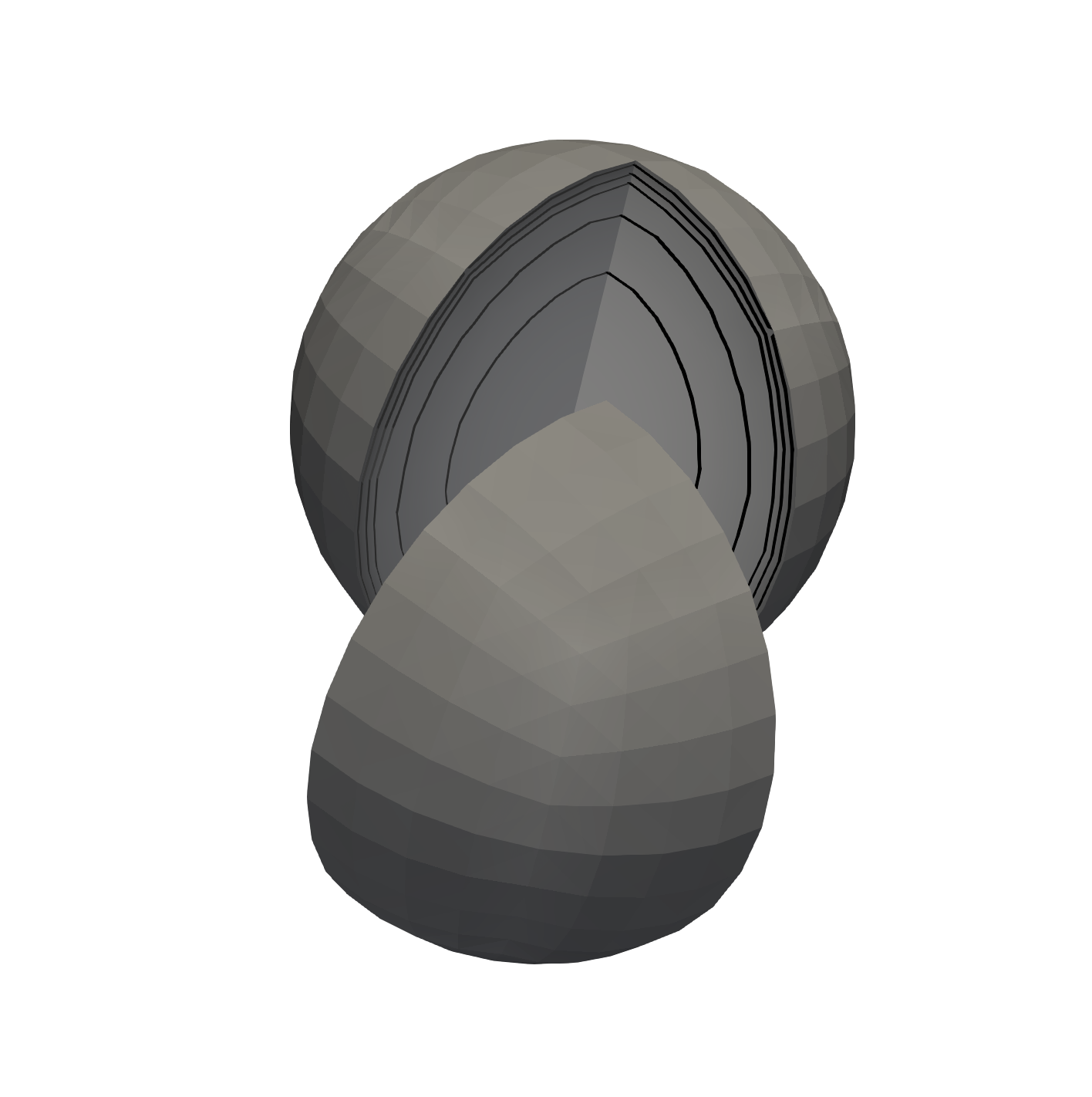}
\caption{Dunford-Integral Method: Numerical approximation of the
  solution to the fractional integral Laplacian for $s=0.3$ and $f=1$ in the unit ball of $\mathbb R^3$ (darker = 0.0, whiter = 0.7).
The lines represent the isosurfaces $\{ u(x) = k/10 \}$ for $k=0,...,7$.
}\label{f:dunf:int:sphere}
\end{figure}

To illustrate the method, we depict in Figure~\ref{f:dunf:int:sphere}
the approximation $U^{k,M}$ for $s=0.3$, $f=1$, and $\Omega = B(0,1)$, the unit ball in $\mathbb R^3$.

%%%%%%%%%%%%%%%%%%%%%%%%%%%%%%%%%%%%%%%%%%%%%%%%%%%%%%%%%%%%%%%%%%%%%%%%%%%%%%%%%%%%
\bibliography{biblio}
\bibliographystyle{plain}
%%%%%%%%%%%%%%%%%%%%%%%%%%%%%%%%%%%%%%%%%%%%%%%%%%%%%%%%%%%%%%%%%%%%%%%%%%%%%%%%%%%%
\end{document}